\documentclass[12]{amsart}

\usepackage{lmodern}  
\usepackage{anyfontsize}

\AtBeginDocument{\fontsize{11}{14}\selectfont}

\usepackage{makecell}

\usepackage{float}

\usepackage{url, fullpage, 
  amssymb,setspace, mathrsfs,fontenc, xcolor}
\usepackage[alphabetic]{amsrefs}
\usepackage{hyperref}
\usepackage{graphicx}
\usepackage[mathcal]{euscript}
\usepackage{verbatim}
\usepackage{cleveref}

\usepackage[inline]{enumitem}

\newtheorem{thm}{Theorem}
\newenvironment{thmbis}[1]
  {\renewcommand{\thethm}{\ref{#1}$'$}
   \addtocounter{thm}{-1}
   \begin{thm}}
  {\end{thm}}
\newtheorem{lem}[thm]{Lemma}

\newtheorem{prop}[thm]{Proposition}

\newtheorem{cor}[thm]{Corollary}
\newtheorem{defe}[thm]{Definition}
 
\theoremstyle{remark}
\newtheorem{rem}[thm]{Remark}

\newtheorem{exam}[thm]{Example}

\let\oldmarginpar\marginpar
\renewcommand\marginpar[1]{\-\oldmarginpar[\raggedleft\footnotesize #1]%
{\raggedright\footnotesize #1}}

\newcommand{\nc}{\newcommand}

\nc{\ssec}{\subsection}

\nc{\on}{\operatorname}

\nc {\cG} {\mathcal{G}}
\nc {\cK}{\mathcal{K}}
\nc {\cD}{\mathcal{D}}
\nc {\cC} {\mathcal{C}}
\nc {\cL} {\mathcal{L}}
\nc {\cE} {\mathcal{E}}
\nc {\cM} {\mathcal{M}}
\nc {\cO}{\mathcal{O}}
\nc {\cF}{\mathcal{F}}
\nc {\cZ}{\mathcal{Z}}

\nc {\bZ}{\mathbb{Z}}
\nc {\bQ}{\mathbb{Q}}

\nc {\uG} {\underline{G}}
\nc {\cB}{\mathcal{B}}
\nc{\rat}{\mathrm{rat}}
      
\nc {\fk}{\mathfrak{K}}
\nc {\fI}{\mathfrak{i}}
\nc {\fg} {\mathfrak{g}}
\nc {\fu} {\mathfrak{u}}
\nc {\fl} {\mathfrak{l}}
\nc {\fn} {\mathfrak{n}}
\nc {\cP} {\mathcal{P}}
\nc {\fz} {\mathfrak{z}}
\nc {\ft} {\mathfrak{t}}
\nc {\fc}{\mathfrak{c}}
\nc {\fh}{\mathfrak{h}}
\nc {\fp}{\mathfrak{p}}

\nc{\tg} {\mathtt{g}}
\nc {\hfg} {\widehat{\fg}}
\nc {\hG} {\widehat{G}}

\nc {\bGm} {\mathbb{G}_m}
\nc{\bC}{\mathbb{C}}
\nc{\bV}{\mathbb{V}}
\nc{\bP}{\mathbb{P}}
\nc{\bA}{\mathbb{A}}

\nc {\mono}{\mathrm{mono}}
\nc{\Sl}{\mathfrak{sl}}

\nc{\ra}{\rightarrow}

\nc {\tU}{\tilde{U}}
\nc {\tSym}{\widetilde{Sym}}

\nc {\Bun}{\mathrm{Bun}}
\nc {\cA}{\mathcal{A}}
\nc {\Fun}{\mathrm{Fun}}
\nc {\crit}{\mathrm{crit}}
\nc {\Ind}{\mathrm{Ind}}
\nc {\Vac}{\mathrm{Vac}}
\nc {\gr}{\mathrm{gr}}
\nc {\ad}{\mathrm{ad}}
\nc {\Sym}{\mathrm{Sym}}
\nc {\Ram}{\mathrm{Ram}}
\nc {\FG}{\mathrm{FG}}
\nc {\Op}{\mathrm{Op}}
\nc {\Hitch}{\mathrm{Hitch}}
\nc {\fb}{\mathfrak{b}}
\nc{\cDt}{\mathcal{D}^\times}
\nc{\cDb}{\mathcal{D}_b^\times}
\nc{\cDbp}{\mathcal{D}_{b'}^\times}
\nc {\gl}{\mathfrak{gl}}
\nc {\Sp}{\mathfrak{sp}}
\nc {\So}{\mathfrak{so}}
\nc {\SP}{\mathrm{Sp}}
\nc {\SO}{\mathrm{SO}}
\nc {\bR}{\mathbb{R}}
\nc {\Hom}{\mathrm{Hom}}

\nc {\Id}{\mathrm{Id}}
\nc {\Tr}{\mathrm{Tr}}
\nc {\rk}{\mathrm{rank}}
\nc {\rank}{\mathrm{rank}}
\nc {\cW}{\mathcal{W}}
\nc {\cI}{\mathcal{I}}
\nc {\Fr}{\mathrm{Fr}}
\nc {\ff}{\mathfrak{f}}

\nc {\LocSys}{\mathrm{LocSys}}
\nc {\Ga}{\mathrm{Ga}}
\nc {\ord}{\mathrm{ord}}
\nc{\pole}{\mathrm{pole}}

\newcommand{\C}{\mathbb{C}}

\newcommand{\Z}{\mathbb{Z}}
\newcommand{\dtt}{\frac{dt}{t}}
\newcommand{\duu}{\frac{du}{u}}

\newcommand{\fs}{\mathfrak{s}}
\newcommand{\Waff}{W^{\mathrm{aff}}}

\newcommand{\om}{\omega}

\DeclareMathOperator{\GL}{GL}
\DeclareMathOperator{\SL}{SL}

\DeclareMathOperator{\FT}{ft}
\DeclareMathOperator{\Res}{\mathrm{Res}}

\DeclareMathOperator{\Lie}{Lie}
\DeclareMathOperator{\Ad}{Ad}
\DeclareMathOperator{\Inn}{Inn}
\DeclareMathOperator{\Out}{Out}
\DeclareMathOperator{\Aut}{Aut}
\DeclareMathOperator{\PAut}{PAut}

\DeclareMathOperator{\rankop}{rank}
\DeclareMathOperator{\diag}{\mathrm{diag}}
\usepackage{chngcntr}
\counterwithout{table}{section}

\newcommand{\sdfrac}[2]{\mbox{\small$\displaystyle\frac{#1}{#2}$}}

\begin{document} 
\title{Differential Galois groups of\\ $G$-connections with Coxeter singularities} 
\author{Masoud Kamgarpour}
\author{Daniel S. Sage}

\begin{abstract}  A fundamental theorem of Katz \cite{Katz87} determines the differential Galois groups of rank $n$ connections on algebraic curves with slope $r/n$ at a singularity, where $\gcd(r,n)=1$. We extend this result to $G$-connections, where $G$ is a simple algebraic group and the slope is $r/h$, with $h$ the Coxeter number of $G$ and $\gcd(r,h)=1$. This allows us to compute the differential Galois groups of a broad class of $G$-connections that have been central to recent advances in the geometric Langlands program and the Deligne--Simpson problem---namely, Coxeter connections, generalised Frenkel--Gross connections, and Airy connections. We apply our results to inverse differential Galois theory by giving uniform and explicit constructions of $G$-connections whose differential Galois groups realise all reductive subgroups of maximal degree.
\end{abstract} 


\address{School of Mathematics and Physics, The University of Queensland, Brisbane, QLD 4072, Australia} 
\email{masoud@uq.edu.au}

\address{Department of Mathematics, University at Buffalo, Buffalo, NY
  14260-2900, USA} \email{dsage@buffalo.edu}

\date{\today}

\maketitle 

\tableofcontents

\section{Introduction} 
\subsection{Background} 
The Galois theory of differential equations, developed by Picard and Vessiot, associates to an $n$th order linear differential equation an algebraic subgroup of $\GL_n$, called the \emph{differential Galois group}.\footnote{Indeed, differential Galois theory was a major impetus for the development of linear algebraic groups in the mid-$20$th century, cf. \cite{BorelHistory}.} This group is finite (resp. solvable) if and only if the solutions are algebraic (resp. Liouvillian). On the other hand, if this group is ``large'', then the solutions are far from being algebraic. Computing differential Galois groups is therefore crucial for understanding special functions such as hypergeometric or Airy functions \cite{Schwarz, Kolchin, Katz87, BH, Katz90, SV, MRS, HMO, BHHW}.

In \cite{Katz87}, N. Katz proved a fundamental result regarding the computation of differential Galois groups. To state this, it is convenient to change perspectives and work with algebraic connections on vector bundles on curves instead of ODEs. 

Let $X$ be a complex smooth projective curve, $S\subset X$ a finite set, and $\nabla$ a rank $n$  connection on $X-S$. The Zariski closure of the monodromy representation $\rho_\nabla:\pi_1(X-S)\to \GL_n$ is known as the \emph{geometric monodromy group} and is denoted by  $G_\nabla^{\mathrm{mono}}$. The differential Galois group $G_\nabla$ contains $G_\nabla^{\mathrm{mono}}$, and they coincide for regular singular connections. For irregular connections, $G_\nabla$ and $G_\nabla^{\mathrm{mono}}$ can be different, as in the case of the exponential connection where $G_\nabla^{\mathrm{mono}}=\{1\}$ and $G_\nabla=\bC^\times$. 

For each singularity of $\nabla$, there is an associated rational number called the \emph{slope}, which is one measure of its irregularity. 

\begin{thm}[\cite{Katz87}]
Let $\nabla$ be a rank $n$ connection on $X-S$ with trivial determinant and connected geometric monodromy. Assume $\nabla$ has an irregular singularity at some $s\in S$ with slope $r/n$, where $\gcd(r,n)=1$. Then $G_\nabla = \mathrm{Sp}_n$ if $n$ is even and $\nabla$ is self-dual; otherwise, $G_\nabla = \SL_n$. \label{t:Katz}
\end{thm}

Katz used this to determine the differential Galois groups of Airy and Kloosterman connections. His results have been crucial in subsequent developments in the field \cite{Katz90, SV, Magid, MS96, CK}.\footnote{We note that Katz also considered connections with nontrivial determinants.} 

Henceforth, let $G$ be a simple complex algebraic group with Lie algebra $\fg$.
Our main goal is to extend Katz's result to the setting of $G$-connections. To this end, we first prove a group-theoretic result which may be of independent interest.

 \subsection{Reductive subgroups of maximal degree}  Recall that the fundamental degrees of $G$ are the degrees of a set of homogeneous algebraically independent generators of $\mathbb{C}[\mathfrak{g}]^G$. The Coxeter number $h = h(G)$ is the largest of these.  More generally, for a connected reductive group $K$, we define the Coxeter number $h(K)$ to be the maximum of the Coxeter numbers of the simple factors of the commutator subgroup $[K,K]$.

\begin{defe} 
A (closed, connected) reductive subgroup $K \subseteq G$ is said to have \emph{maximal degree} if $h(K)=h(G)$. 
\end{defe}

\begin{thm} \label{t:largeIntro}
 Let $G$ be a simple group. Then, reductive subgroups of maximal degree in $G$ are simple and arise as  fixed points of pinned automorphisms of $G$. 
\end{thm} 

For a more precise statement of this result, as well as a table listing the root systems of reductive subgroups of maximal degree, we refer the reader to \S \ref{s:large}.

\subsection{Main result} The notions of monodromy, slope, and differential Galois group extend naturally from the setting of connections on vector bundles to that of connections on $G$-bundles on curves, cf. \cite{FG, CK, BSminimal}. By definition, the differential Galois group  of a $G$-connection is a closed subgroup of $G$. 

\begin{defe} \label{d:Cox}
A $G$-connection on $X-S$ is said to have a \emph{Coxeter singularity} at $s\in S$ if it has an irregular singularity at $s$ with slope $r/h$, where $\gcd(r, h) = 1$. 
\end{defe}

As discussed below, Coxeter singularities have a highly rigid structure. Our main result shows that the presence of a single Coxeter singularity severely restricts the differential Galois group. 

Let $\Aut(G)$ denote the group of (algebraic) automorphisms of $G$, $\Inn(G)$ the subgroup of inner automorphisms, and $\Out(G):=\Aut(G)/\Inn(G)$. Recall that $\Out(G)$ is a finite group isomorphic to $\{1\}$, $\Z/2\Z$, or $S_3$. 
For each $\sigma\in \Aut(G)$, we let $\overline{\sigma}$ denote its image in $\Out(G)$.

\begin{thm} \label{t:mainIntro} Let  $G$ be a simple group and $\nabla$ a $G$-connection on $X-S$ with a Coxeter singularity at some $s\in S$. Assume the geometric monodromy  group $G_\nabla^{\mathrm{mono}}$ is connected. Then the differential Galois group $G_\nabla\subseteq G$ is a reductive subgroup of maximal degree whose type is given in Table \eqref{ta:mainIntro}.  
\begin{table}[h]
\centering
\renewcommand{\arraystretch}{1.5} 
\begin{tabular}{|l|l|} 
\hline
$G$  &  $G_\nabla$ \\ \hline
$A_1$ & $A_1$\\\hline 
$A_{2n-1},\ n \geq 2$ & 
\begin{tabular}[t]{@{}l@{}}
$C_{n}$ if $\sigma(\nabla)\simeq \nabla$ for some $\sigma\in \Aut(G)$ with $\overline{\sigma}\neq 1$. \\
$A_{2n-1}$ otherwise
\end{tabular} \\ \hline
$A_{2n}, \ n\geq 1$ & $A_{2n}$ \\\hline
$B_3 $ & 
\begin{tabular}[t]{@{}l@{}}
$G_2$ if $\tau(\nabla)\simeq \nabla$ for some $\tau\in \Aut(G')$ with $\overline{\tau}$ of order $3$ \\
$B_3$ otherwise
\end{tabular} \\ \hline
$B_n,\ n\geq 4$ & $B_n$\\\hline 
$C_n,\ n\geq 2$ & $C_n$\\\hline 
$D_4$ & 
\begin{tabular}[t]{@{}l@{}@{}}
$G_{2}$ if $\tau(\nabla)\simeq \nabla$ for some $\tau\in \Aut(G)$ with $\overline{\tau}$ of order $3$\\
$B_3$ if $\#\{ \overline{\sigma} \, | \, \sigma(\nabla)\simeq \nabla,\,\, \sigma\in \Aut(G) \}=2$\\
$D_4$ otherwise
\end{tabular} \\ \hline
$D_{n},\ n \geq 5$ & 
\begin{tabular}[t]{@{}l@{}}
$B_{n-1}$ if  $\sigma(\nabla)\simeq \nabla$ for some $\sigma\in \Aut(G)$ with $\overline{\sigma}\neq 1$\\
$D_n$ otherwise
\end{tabular} \\ \hline
$E_{6}$ & 
\begin{tabular}[t]{@{}l@{}}
$F_4$ if  $\sigma(\nabla)\simeq \nabla$ for some $\sigma\in \Aut(G)$ with  $\overline{\sigma}\neq 1$\\
$E_6$ otherwise
\end{tabular} \\ \hline
$E_7$ & $E_7$\\\hline 
$E_8$ & $E_8$\\\hline 
\end{tabular}
\vspace{0.8em} \caption{Differential Galois groups of connections with a
  Coxeter singularity and connected geometric monodromy.  In type $B_3$, $G'$  is taken to be a group of type $D_4$ containing $G$, and we think of $\nabla$ as a $G'$-connection via the embedding $G\hookrightarrow G'$.}
\label{ta:mainIntro} 
\end{table}
\end{thm} 

\begin{rem} \label{r:forAll} 
As explained in \S \ref{ss:Auto}, one can replace every instance of ``for some" by ``for all" in Table \ref{ta:mainIntro}.  In particular, when $G=\SL_{2n}$ and $\sigma(A)=(A^t)^{-1}$, we recover Theorem \ref{t:Katz}. 
\end{rem} 

 In general, deciding whether $\sigma(\nabla)$ is isomorphic to $\nabla$ can be challenging. Nevertheless, as discussed below, this can be done for several important classes of connections, including generalised Frenkel--Gross and Airy connections.

\subsection{Normal forms of Coxeter singularities} Let $\cDt=\mathrm{Spec}(\bC(\!(t)\!))$ denote the formal punctured disk. 

\begin{defe}\label{d:formCox} A formal Coxeter $G$-connection is a $G$-connection on $\cD^\times$ with an irregular singularity of slope $r/h$, where $r$ is a positive integer satisfying $\gcd(r, h) = 1$.  
\end{defe}

The first step in the proof of Theorem \ref{t:mainIntro} is to compute the Jordan canonical form (also known as the Levelt--Turrittin normal form) of a formal Coxeter connection. This allows us to determine its (local) differential Galois group. A key consequence is that this group is irreducible. Next, by a theorem of Gabber (see  \cite{Katz87}), the fact that $G_\nabla^\mathrm{mono}$ is connected implies that so is $G_\nabla$. Since $G_\nabla$ contains the local differential Galois group, it follows that $G_\nabla$ is irreducible. By a fundamental result of Serre and Borel--Tits, $G_\nabla$ is then reductive. Viewing $\nabla$ as a $G_\nabla$-connection, we conclude that one of the fundamental degrees of $G_\nabla$ must be equal to $h$. This, in turn, implies that $G_\nabla$ is a reductive subgroup of maximal degree.

To construct and analyse explicit examples of connections to which our theorem applies, we use the \emph{rational canonical form} of formal Coxeter connections. Let $S$ be a maximal torus in the loop group of $G$. A formal connection is called \emph{$S$-toral} if the leading term of its matrix, with respect to some Moy--Prasad filtration, is regular semisimple with centraliser conjugate to $S$. The key property of $S$-toral connections is that they can be ``diagonalised'' into the Cartan subalgebra $\fs = \Lie(S)$. A maximal torus is called a \emph{Coxeter torus} if it corresponds to the Coxeter class in the Weyl group, under the natural bijection between conjugacy classes in the Weyl group and conjugacy classes of maximal tori in the loop group. We will show that a formal Coxeter connection is Coxeter toral. This implies that such a connection admits an especially simple rational canonical form.\footnote{After the first draft of this paper appeared on the arXiv in 2023, we were informed that some of these local results were obtained independently in \cite{JY}, \cite{ChenLi}, and  \cite{Lingfei}.}

\subsection{Applications to Coxeter connections} 

\begin{defe} \label{d:CoxConn}
A \emph{Coxeter $G$-connection} on $\bGm$ is a $G$-connection with a regular singularity at $\infty$ and an irregular singularity at $0$ of slope $r/h$, where $r$ is a positive integer satisfying $\gcd(r,h)=1$.
\end{defe}

Coxeter connections were introduced in \cite{KSCoxeter} and have played a central role
in recent advances on the Deligne--Simpson problem \cite{KLMNS, JY, Sagelocal} and the
wildly ramified geometric Langlands program \cite{KML, ChenLi, YiAiry}.\footnote{In
the literature, the regular singularity is usually placed at $0$ and the irregular one at $\infty$.}

We will use Theorem \ref{t:mainIntro} to compute the differential Galois groups of a wide class of Coxeter connections. 
Here, we restrict ourselves to the particularly important cases of generalised Frenkel--Gross connections and Airy
connections.  

\begin{defe} \label{d:FG}
A generalised Frenkel--Gross connection is a (Coxeter) connection on $\bGm$ with an irregular singularity at $0$ of slope $1/h$ and 
 a regular singularity at $\infty$. 
\end{defe} 

Let us fix a pinning $\cP=(B, T, \{x_\alpha\})$ for $G$, and let $N$ be the corresponding regular nilpotent element defined in \S \ref{ss:Notation}. We choose in addition a generator $E$ of the highest root space. 
As we shall see in \S\ref{ss:FG}, every generalised Frenkel--Gross connection can be written in the form
\begin{equation}\label{eq:FG}
\nabla = d+ a(N + S + t^{-1}E)\,\dtt,
\end{equation} 
where $S$ and $a$ are arbitrary elements of $\ft:=\Lie(T)$ and $\bC^\times$, respectively.

Let us write $\Aut(G, \cP)$ for the subgroup of $\Aut(G)$ consisting of automorphisms which preserve $\cP$. Recall that $\Aut(G, \cP)$ is isomorphic to  $\Out(G)$. (It follows that in Table \ref{ta:mainIntro}, we can restrict to $\cP$-pinned automorphisms.)  Let $\widetilde{W}=W\ltimes X_*(T)$ denote the extended affine Weyl group of $G$.

\begin{thm}\label{t:FG} Let $S\in \ft$ be a torsion-free element as in Definition \ref{d:torsionFree}. Then
the differential Galois group of the generalised Frenkel--Gross connection \eqref{eq:FG} is given by a simplified version of Table~\eqref{ta:mainIntro}, in which the conditions $\sigma(\nabla) \simeq \nabla$ are replaced by restricting $\sigma$ to $\Aut(G, \cP)$ and requiring that $\sigma(S)$ be $\widetilde{W}$-conjugate to $S$ (and similarly for $\tau$ in types $B_3$ and $D_4$).
\end{thm}

The significance of this theorem is that it reduces the difficult problem of determining whether $\nabla$ is isomorphic to $\sigma(\nabla)$ to the much simpler question of whether $S$ is $\widetilde{W}$-conjugate to $\sigma(S)$. In particular, the latter can readily be verified computationally. Note that when $S = 0$, we recover the differential Galois group of the original Frenkel--Gross connection \cite{FG}. However, our argument does not use the fact that the monodromy is regular, which was a crucial part of the original computation.

\begin{defe} \label{d:Airy} An Airy $G$-connection is a (Coxeter) connection on $\bP^1-\{0\}\simeq \bA^1$ with an irregular singularity at $0$ of slope $1+1/h$. 
\end{defe} 
As we shall see in \S\ref{ss:Airy}, every Airy $G$-connection can be written in the form
\begin{equation} \label{eq:Airy}
\nabla = d + at^{-1}(N + X + t^{-1}E)\,\dtt, 
\end{equation}
where $N$, $E$, and $a$ are as above, and $X$ is an element of the Borel subalgebra $\fb:=\Lie(B)\subset \fg$. 

Let $\nabla_0$ denote  the formal Coxeter connection on the formal disk $\cDt_0$ around zero obtained by restricting $\nabla$.  Let $\FT(\nabla_0)$ be the formal type of $\nabla_0$ as defined in \S \ref{ss:Rational}.\footnote{Here, we normalise the formal type to have leading term $at^{-1}(N+t^{-1}E)$.}  (Note that $d+\FT(\nabla_0)\dtt$ is the rational canonical form of $\nabla_0$.)

\begin{thm}\label{t:Airy}
The differential Galois group of the Airy $G$-connection \eqref{eq:Airy} is given by a simplified version of Table~\eqref{ta:mainIntro}, in which the conditions $\sigma(\nabla) \simeq \nabla$ are  replaced by restricting $\sigma$ to $\Aut(G, \cP)$ and requiring  $\FT(\sigma(\nabla_0))=\FT(\nabla_0)$ (and similarly for $\tau$ in types $B_3$ and $D_4$).
\end{thm}

The significance of this theorem is that it reduces a global condition to a local one. In \S \ref{s:algo}, we present an algorithm for determining the rational canonical form of a formal Coxeter connection. This algorithm can be implemented in practice, and thus, the conditions of the theorem may be  verified explicitly.

\subsection{Applications to inverse differential Galois theory} It is a classical result, going back to Hilbert, that the Galois group of a generic degree $n$ polynomial with rational coefficients is  $S_n$. Refinements of this result continue to attract considerable interest; see, for example, \cite{Bhargava}. Explicit examples valid for all $n$, however, remain scarce. Perhaps the best-known example is the polynomial $x^n - x - 1$, which Selmer proved to be irreducible for all $n$ \cite{Selmer}, and which was later shown to have Galois group $S_n$ \cite{NV}. It is not easy to generalise this result, even to other trinomials; see, for instance, \cite{Osa87a,Osa87b,CMS97,CMS99}, where explicit criteria are given for an irreducible degree $n$ trinomial to have a ``large'' Galois group, namely one containing the alternating group $A_n$.

Analogous questions arise in differential Galois theory. As in the
arithmetic case, it is easier to show the existence of $G$-connections
with differential Galois group $G$ than to exhibit explicit examples
of them.  Mitschi and Singer \cite{MS96} gave a uniform existence
result, constructing linear differential equations with differential
Galois group $G$ for every connected semisimple group $G$, using
 ``Chevalley modules". These are faithful $G$-modules
with no $G$-stable lines, but such that every proper connected
subgroup stabilises a line. These modules
are difficult to describe explicitly, especially in higher rank, and
no uniform explicit construction is currently available. As a result,
this approach does not typically produce concrete $G$-connections. By contrast, in \S\ref{s:applications}, we give a uniform and fully
explicit construction of Coxeter $G$-connections realising all the possibilities listed in Table \eqref{ta:main}, including $G$ itself. 
 More precisely, we will prove: 

\begin{thm}\label{t:Galois}
For every positive integer $r$ coprime to $h$ and every reductive subgroup $K \subseteq G$ of maximal degree, there exists an explicitly and uniformly defined family of Coxeter $G$-connections on $\bGm$ with slope $r/h$ and differential Galois group $G_\nabla = K$. If $r > h$, such a family exists on $\bA^1$.
\end{thm}

This theorem is new even for $r=1$, since  the differential Galois groups of generalised Frenkel–Gross connections had not previously been computed.

\subsection{Relevance in the geometric Langlands program}
Let $\hG$ denote the Langlands dual group of $G$. The central conjecture of the geometric Langlands program (cf.~\cite{BZN}) asserts that every $G$-connection $\nabla$ on $X - S$ arises as the eigenvalue of a Hecke eigensheaf on $\mathrm{Bun}_{\hG}(X)$, with additional structure  at the singular points $S$. 

If the differential Galois group $G_\nabla\subseteq G$ is a proper connected reductive subgroup with Langlands dual group $\widehat{G_\nabla}$, then the conjecture implies that there also exists a Hecke eigensheaf for $\nabla$ on $\mathrm{Bun}_{\widehat{G_\nabla}}(X)$. This relationship between Hecke eigensheaves on $\mathrm{Bun}_{\hG}(X)$ and $\mathrm{Bun}_{\widehat{G_\nabla}}(X)$ is known as \emph{functoriality}. On the other hand, the Hecke eigensheaves associated to $\nabla$'s with $G_\nabla=G$ do not come from smaller groups and are expected to form the basic building blocks of the constructible derived category of $\mathrm{Bun}_{\hG}$.

In the irregular setting, Hecke eigensheaves are currently known in only a small number of  cases, most notably for certain Coxeter $G$-connections such as generalised Frenkel--Gross and Airy connections. This underscores the importance of determining their differential Galois groups.

\subsection{Relevance in the Deligne--Simpson problem}  For each $s\in S$, fix a (formal) $G$-connection
$\nabla_s$ on the formal punctured disk $\cDt_s$ around $s$. Roughly
speaking, the Deligne--Simpson problem asks whether 
there exists an irreducible $G$-connection $\nabla$ on $X-S$ such that for each $s\in S$, the restriction of $\nabla$ to $\cDt_s$ is isomorphic to $\nabla_s$. (For a more precise formulation, see \cite{KLMNS, Sagelocal}.) 
Restricting to Fuchsian connections on $X=\bP^1$, this problem takes on
an especially simple form: For each $s\in S$, fix a conjugacy class
$\cC_s\subset G$. Does there exist an irreducible collection of $X_s \in \cC_s$ satisfying 
\[
\prod_{s\in S} X_s=1? 
\]

When $G=\GL_n$ and $X=\bP^1$, Kostov formulated an additive version
of this problem~\cite{Kostov03}; it was solved by Crawley-Boevey
\cite{CB} by building a bridge to quiver varieties.   Hiroe also used
quiver varieties in extending Crawley-Boevey's results to include
``unramified'' irregular singular points~\cite{Hiroe}.  
Much less is known when ramified irregular singular points are allowed
and $G$ is arbitrary.  As discussed in \cite{KLMNS, LSN, JY, Sagelocal}, the ramified Deligne–Simpson problem is best understood for the class of $G$-connections that are central to this paper, namely Coxeter connections.

\subsection{Future directions} 
Let us call a $G$-connection on $\bGm$ \emph{elliptic} if it has an irregular singularity of slope $r/m$, where $\gcd(r, m) = 1$ and $m$ is a regular elliptic number for $G$. Most aspects of our theory extend naturally to the setting of elliptic $G$-connections. In particular, if $\nabla$ is an elliptic connection with connected geometric monodromy, then $G_\nabla$ is a connected reductive subgroup of $G$, one of whose fundamental degrees divides $m$. A remaining challenge is to classify all such subgroups. We plan to return to this question in future work.

In a different direction, let $\nabla$ be an irreducible rank $n$ connection on a smooth curve. Suppose that $\nabla$ has an irregular singularity with slope equal to $a/b$, where  $\gcd(a, b) = 1$ and that this slope appears exactly $b$ times. Then,  Katz has proved that $G_\nabla$ is either $\mathrm{SO}_n$, $\mathrm{Sp}_n$, or $\mathrm{SL}_n$, with a few exceptions when $n = 7$, $8$, or $9$ \cite[Theorem 2.8.1]{Katz90}. It would be interesting to find an analogue of this result for $G$-connections.

\subsection{Structure of the paper}

In \S\ref{s:group}, we collect results on reductive groups. In \S \ref{s:large}, we study reductive subgroups of maximal degree and prove (a sharper version of)  Theorem \ref{t:largeIntro}.  In \S\ref{s:connections}, we describe the Jordan canonical form of formal Coxeter $G$-connections and compute their  differential Galois groups. In \S\ref{s:global}, we prove Theorem \ref{t:mainIntro}.  The rational canonical forms and moduli spaces of formal Coxeter connections are described in \S \ref{s:rational}.  In \S \ref{s:applications}, we determine the differential Galois groups of Coxeter connections and prove Theorems \ref{t:FG}, \ref{t:Airy}, and \ref{t:Galois}. 

\subsection{Acknowledgments} MK was supported by the Australian Research Council, and DS was supported by the National Science Foundation and the Simons Foundation. We would like to thank Konstantin Jakob and Lingfei Yi for many helpful conversations and their interest in this project. We also thank M. van der Put for answering our questions about generic differential Galois groups.


\section{Group-theoretic preliminaries} \label{s:group} 
In this section, we collect several foundational results on reductive groups that will be used throughout the paper. In \S\ref{ss:pinned} (resp. \S \ref{ssAuto}) we recall  some  facts about pinned (resp. semisimple) automorphisms of simple groups. In \S \ref{sss:PinnBig}, we prove a result regarding a characterisation of pinned automorphisms. 
In \S\ref{ss:connected}, we establish a criterion for when the Zariski closure of a cyclic subgroup is connected. Next, in \S\ref{ss:irreducible}, we show that any connected subgroup containing a lift of a Coxeter element is necessarily reductive. In the subsequent \S\ref{ss:CoxeterGrading}, we study the Coxeter grading. Finally, in \S\ref{ss:loop} we review basic material on loop algebras, and in \S\ref{ss:Weyl} we recall some facts about affine Weyl groups.

\subsection{Notation} \label{ss:Notation} 
 Throughout the paper, $G$ denotes a simple complex algebraic group, i.e., a connected complex algebraic group whose Lie algebra $\mathfrak{g}$ is simple.  

\subsubsection{} We fix a pinning for $G$, i.e., a triple 
\[
\mathcal{P}:=(B, T, \{x_\alpha\}), 
\]
 where $B$ is a Borel subgroup, $T$ is a maximal torus
in $B$, and, for each simple root $\alpha$, $x_\alpha$ is a nonzero
vector in the corresponding root space. We let $\Phi$ and $\Delta$ denote the
set of roots and simple roots, respectively. 
Let $\fb$ and $\ft$ denote the Lie algebras of $B$ and $T$
 and  $W$ the Weyl group. 

\subsubsection{} 
Using the standard $\mathfrak{sl}_2$-triple associated to each simple root,
our chosen pinning determines a nonzero element $y_\alpha \in \fg_\alpha$ for every negative simple root $\alpha$,

We extend this choice to an \emph{affine pinning} by selecting a nonzero element
$E$ in the highest root space.
We then define
\[
N := \sum_{\alpha \in -\Delta} y_\alpha .
\]
By a theorem of Kostant, $N$ is a regular nilpotent element of $\fg$,
and the sum $N+E$ is regular semisimple.

\subsubsection{} For each $g\in G$, we let $c_g: G\ra G$ denote the inner automorphism 
\[
c_g(x)=gxg^{-1},\qquad x\in G. 
\]

\subsubsection{} 
We note that many of our constructions and results apply (after mild modifications) if $G$ is assumed to be a connected semisimple group over an algebraically closed field of characteristic
$0$. In particular, this is the case up through \S \ref{sss:NGT}.

\subsection{Pinned automorphisms}  \label{ss:pinned}
A reference for the material in this section is \cite[\S 23]{Milne}. Let $\Aut(G)$ denote the group of  automorphisms of $G$ and $\Inn(G)$ the subgroup of inner automorphisms. Then $\Inn(G)\simeq G^{\mathrm{ad}}$, where $G^{\mathrm{ad}}:=G/Z(G)$ is the corresponding adjoint group.  We have a canonical exact sequence 
\begin{equation}\label{eq:exact}
1\ra \Inn(G) \ra \Aut(G) \ra \Out(G)\ra 1. 
\end{equation}

\subsubsection{$\cP$-pinned automorphisms} 
 Let $\Aut(G,\cP)\subseteq \Aut(G)$ denote the subgroup of
 \emph{$\cP$-pinned automorphisms} of $G$, i.e. automorphisms of $G$
 which preserve $\cP$.  Explicitly, this means that the automorphism
 preserves $B$ and $T$ and permutes the $x_\alpha$'s. The composition 
 \[
 \Aut(G,\cP) \hookrightarrow \Aut(G) \twoheadrightarrow \Out(G) 
 \]
 is an isomorphism. Thus, $\Aut(G, \cP)$ provides a section of the exact sequence \eqref{eq:exact} and we have: 
 \[
 \Aut(G) \simeq  \Inn(G) \rtimes \Aut(G,\cP)\simeq G^{\mathrm{ad}} \rtimes \Aut(G,\cP). 
 \]

\subsubsection{} 
If $\cP'=(B', T', \{x'_\alpha\}), $ is another pinning, then there
exists $g\in G$ such that $c_g(B)=B'$, $c_g(T)=T'$, and
$\Ad_g(x_\alpha)=x'_\alpha$ for all $\alpha\in \Delta$.  Thus, pinned automorphisms with respect to $\cP$ and $\cP'$ are conjugate; in other words, the map
\[
\Aut(G,\cP) \ra \Aut(G, \cP'),\qquad \sigma \mapsto c_g \circ \sigma \circ c_{g^{-1}},
\] 
is an isomorphism. Indeed, it is immediate that $c_g \circ \sigma
\circ c_{g^{-1}}$ preserves $T'$ and $B'$ while 
\[
\Ad_g\circ \, \sigma \circ \Ad_{g^{-1}}(x'_\alpha)=\Ad_g(\sigma(x_\alpha))=\Ad_g(x_{\sigma(\alpha)})=x'_{\sigma(\alpha)}.
\]

\subsubsection{} 
By a pinned automorphism of $G$, we mean an automorphism which fixes some, but not necessarily our chosen, pinning. Let $\PAut(G)\subseteq \Aut(G)$ denote the subset consisting of automorphisms that are pinned with respect to some pinning.
 By the above discussion
\[
\PAut(G)=\bigcup_{g\in G} c_g \Big(\Aut(G, \cP)\Big) c_{g^{-1}}.
\]

\subsubsection{} 
Let $\Gamma_G$ denote the group of automorphisms of the Dynkin diagram of $G$. Then we have an isomorphism 
\[
 \Out(G)\simeq \Gamma_G.
\]
It is known that $\Gamma_G$ is isomorphic to $S_3$ in type $D_4$, isomorphic to $\mathbb{Z}/2$ in types $A_n$ for $n\geq 2$,  $D_n$ for $n \geq 5$, and $E_6$, and is trivial otherwise.

\subsubsection{} \label{sss:fixedPinn} Let $\sigma\in \PAut(G)-\{1\}$. 
The possibilities for the types of the fixed point
subgroups $(G^\sigma)^\circ \subseteq G$ are listed in Table \ref{ta:pinned}. 
By Dynkin's classification, these subgroups are maximal except for $G_2\hookrightarrow D_4$.\footnote{See \S \ref{sss:maximal} for a direct proof of this fact in all types except $A_{2n}$.} In addition, in all cases except $B_3\hookrightarrow D_4$, the subgroups in Table \ref{ta:pinned} are unique up to $G$-conjugacy. 
By triality, there are three conjugacy classes of embeddings $B_3\hookrightarrow D_4$ which are related by outer automorphisms.

\begin{table}[h]
\centering
\renewcommand{\arraystretch}{1.5}
\[
\begin{array}{|c|c|}
\hline
G & (G^\sigma)^\circ  \\ \hline
A_{2n}, \, n\geq 1  & B_n  \\ \hline
A_{2n-1}, \, n\geq 2 & C_n    \\ \hline
D_{n}, \, n\geq 5  & B_{n-1}       \\ \hline
D_4      & G_2\quad   \text{if } \sigma^3=1 \\ \hline
D_4      & B_3  \quad  \text{if } \sigma^2=1 \\ \hline
E_6      & F_4      \\ \hline
\end{array}
\]
\caption{Fixed-point subgroups under nontrivial pinned automorphisms}
\label{ta:pinned}
\end{table}

\subsubsection{} The following proposition shows that the subgroups in the above table can only arise as fixed points of \emph{pinned} automorphisms. 
\begin{prop}\label{p:PinnBig} If $\tau\in\Aut(G)$ has the property that $\fg^\tau\simeq \fg^\sigma$ for some pinned automorphism $\sigma$, then
  $\tau$ is itself a pinned automorphism.
\end{prop}

Note that this proposition gives another characterisation of the subset $\PAut(G)\subseteq \Aut(G)$. We give the proof in \S \ref{sss:PinnBig}.

\subsection{Semisimple automorphisms} \label{ssAuto}
Recall that we have a Jordan decomposition for elements of the algebraic group $\Aut(G)$; thus, we can speak of {semisimple} and unipotent automorphisms. Following \cite{SteinbergEnd}, we call $\tau\in \Aut(G)$ \emph{quasi-semisimple} if $\tau$ fixes a pair $(B', T')$ consisting of a Borel subgroup $B'$ containing a maximal torus $T'$.  Note that every pinned automorphism is quasi-semisimple, but the converse is not true. We now describe how these notions are related.

\begin{thm}\label{t:ssAuto} The following are equivalent: 
\begin{enumerate} 
\item[(i)] $\tau$ is semisimple.
\item[(ii)] $\tau$ is quasi-semisimple.
\item[(iii)] $\tau$
  is the composition of a $\cP'$-pinned automorphism and conjugation by an
  element of $T'$, for some pinning $\cP'=(B', T', \{x_\alpha'\})$.
  \item[(iv)] $(G^\tau)^\circ$ is reductive. 
\end{enumerate} 
\end{thm}

\begin{proof} 
$(i) \iff(ii)$:  This is established in \cite[\S
 9]{SteinbergEnd}. (In positive characteristic, we only have (i) implies (ii).)

$(ii) \implies (iii)$: Suppose $\tau$ is quasi-semisimple and fixes the pair
 $(B', T')$.   Then $\tau$ induces an automorphism
 $\overline{\tau}$ of the
 Dynkin diagram, and there exists a (unique) element 
 $\pi\in\Aut(G,\cP')$ such that $\overline{\pi}=\overline{\tau}$.  By a theorem of Kostant, we can choose
 $t\in T'$ such that
 \[
 \Ad_t(x'_{\overline{\tau}(\alpha)})=\tau(x_\alpha')
 \]
  for all simple roots $\alpha$.  Now, the
 automorphisms induced on $\fg$ by $c_t\circ\pi$ and $\tau$ agree on
 $\ft'$ and the simple roots spaces, hence they coincide on $\fg$.  It
 follows that $\tau=c_t\circ\pi$, implying (iii). 
 
 $(iii) \implies (ii)$: This is obvious because $\tau$ preserves $(B', T')$.

$(ii)\implies (iv)$: This is proved in \cite[\S
8]{Steinberg} and \cite{DM}. 

$ (iv) \implies (i)$: Suppose 
$G':=(G^\tau)^\circ$ is reductive. Our aim is to show that $\tau$ is semisimple.

Without loss of generality, we may assume $G$ is of adjoint type, i.e., the centre of $G$ 
is trivial. Let $\tau=\sigma \upsilon$ be the Jordan decomposition in
$\Aut(G)$, where $\sigma$ is semisimple and $\upsilon$ is
unipotent. By  Lemma \ref{l:unipotent} below, there exists a unipotent
element $u\in G'$ such that $\upsilon$ is conjugation by $u$. We will show that $u=1$.

Suppose $u\neq 1$. 
By a theorem of Borel and Tits \cite{BorelTits}, there exists a parabolic subgroup $P'\subset G'$
such that $u$ belongs to the unipotent radical $U'\subset P'$ and
$(G')^u\subset P'$. Now, $U'\cap (G')^u$ is a nontrivial unipotent (connected) normal subgroup of $(G')^u$.
Hence, $\Lie((G')^u)=(\fg')^\upsilon=(\fg^\tau)^\upsilon$ is not reductive. However, we will show  that $(\fg')^\upsilon = \fg^\tau$ and hence, is reductive. This is a contradiction.

Finally, observe that
$\Aut(G)=\Aut(\fg)$ embeds naturally in $\GL(\fg)$.  Set $V=\fg^\tau$.  Since
$\tau$ preserves $V$, so do $\sigma$ and $\upsilon$, and moreover,
$\tau\vert_V=\sigma\vert_V\upsilon\vert_V$ is the Jordan composition
of $\tau\vert_V$ \cite[\S I.4, Corollary 1]{Borel}.  Since $\tau\vert_V=1_V$, it follows that
$\sigma\vert_V=1_V=\upsilon\vert_V$.  Thus,
$\fg^\tau\subseteq\fg^\sigma\cap\fg^\upsilon$, so $(\fg^\tau)^\upsilon = \fg^\tau$. 

\end{proof} 

 \begin{lem}\label{l:unipotent} 
  Let $G$ be a group of adjoint type. Let $\tau$ be an automorphism of $G$ with Jordan decomposition $\tau=\sigma\upsilon$, where $\sigma$ is semisimple and $\upsilon$ is unipotent.  Then there exists $u\in (G^\tau)^\circ$ such that $\upsilon$ is conjugation by $u$. \end{lem} 
 
\begin{proof} Since every element of the finite group
$\Out(G)$ is semisimple, we must have $\overline{\upsilon}=1$, where
$\overline{\upsilon}$ is the image of $\upsilon$ under the map $\Aut(G)\ra
\Out(G)$. In other words, 
\[
\upsilon\in\Aut(G)^\circ =\Inn(G)\simeq G.
\]
Thus, $\upsilon=c_u$ for some unipotent element
$u\in G$.  It follows that the automorphism $\tau\in \Aut(G)$ is  given by 
\[
\tau(g) = \sigma(ugu^{-1}),\qquad \forall \, g\in G. 
\]
As $\sigma$ and $c_u$ commute, we have 
\[
\sigma(u)\sigma(g) \sigma(u)^{-1} = \sigma(ugu^{-1}) = \sigma c_u(g) = c_u \sigma(g) = u \sigma(g) u^{-1}. 
\]
Thus, $u^{-1}\sigma(u) \sigma(g) \sigma(u)^{-1} u= \sigma(g)$. Since this holds for all $g\in G$, we conclude that 
\[
u^{-1}\sigma(u)\in Z(G)=\{1\}.
\] 
Thus, $\sigma(u)=u$, i.e., $u$ belongs to the fixed points $G^{\sigma}\cap G^\upsilon\subseteq G^\tau$. As $u$ is unipotent, this implies that $u$ is in the neutral component  $(G^{\tau})^\circ$. 
\end{proof}

\subsection{Proof of Proposition
  \ref{p:PinnBig}} \label{sss:PinnBig}  
  We want to show that if $\fg^\tau\simeq \fg^\sigma$, and $\sigma$ is pinned, then
  so is $\tau$.  The case $\sigma=1$ is
trivial, so we assume $\sigma\neq 1$; thus, $\fg^\tau$ is a proper
subalgebra. 

\subsubsection{} As noted in \S \ref{sss:fixedPinn}, every subalgebra of $\fg$ isomorphic to $\fg^\sigma$ is conjugate
to $\fg^{\sigma'}$ for some pinned automorphism $\sigma'$. (Indeed, in all cases except $B_3\hookrightarrow D_4$, one
can take $\sigma'=\sigma$.)  Thus, we may
assume without loss of generality that $\fg^\tau=\fg^\sigma$.  We will prove that $\tau=\sigma$ in all types, except when $\fg$ has type $D_4$ and $\sigma$ has order three, in which case, we have $\tau=\sigma$ or $\tau=\sigma^2$. 

\subsubsection{} We first show that $\tau$ cannot be inner. 
As $\fg^\sigma$ is simple, Theorem~\ref{t:ssAuto} implies that $\tau$ is
  semisimple. If $\tau$ were inner, then it would have to be conjugation by a semisimple element of $G$, and $\fg^\tau$ would have maximal rank. This contradicts the fact that $\rankop(\fg^\sigma)<\rankop(\fg)$; see Table \ref{ta:pinned}.

\subsubsection{} 
Suppose we are not in type $D_4$. Then,  $\overline{\tau} =
  \overline{\sigma}$. Thus, $\pi=\sigma^{-1}\tau$ is
 an inner automorphism, and  $\fg^{\pi}$ contains
 $\fg^\tau=\fg^\sigma$.  As noted in \S \ref{sss:fixedPinn}, the subalgebra $\fg^\sigma$ is maximal, so
 either $\fg^\pi=\fg$ or $\fg^\pi=\fg^\sigma$. The latter case is not
 possible.  Indeed, Theorem~\ref{t:ssAuto} shows that $\pi$ is an inner
 semisimple automorphism.  This implies that $\rankop(\fg^\pi)=\rankop(\fg)>\rankop(\fg^\sigma)$, which is a contradiction.  It follows that $\fg^\pi=\fg$, so
 $\pi=1$, and $\tau=\sigma$. 
 
 \subsubsection{} Henceforth, we consider type $D_4$. We first show that
$\overline{\tau}$ and $\overline{\sigma}$ have the same order.     By a
  theorem of Kac~\cite{Kac}, the rank of $\fg^\tau$
equals the number of orbits of the action of $\overline{\tau}\in
\Out(G)=\Gamma_G$ on the set of simple roots of $\fg$. Thus, the rank
of $\fg^\tau$  equals $2$ (resp. $3$) if $\overline{\tau}$ has order
$3$ (resp. $2$) and similarly for $\sigma$. As $\fg^\tau=\fg^\sigma$, we conclude that $\overline{\tau}$ and $\overline{\sigma}$ have the same order.

\subsubsection{} 
We claim that if $\overline{\tau}$ and $\overline{\sigma}$ have order $2$, then 
$\overline{\tau}=\overline{\sigma}$.  If not, choose a pinned
involution $\sigma'$ with respect to the same pinning as $\sigma$ for which
$\overline{\tau}=\overline{\sigma'}$.  Note that $\rho=\tau^{-1}\sigma'$ is
an inner automorphism, and $\fg^\rho$ contains
$\fg^\sigma\cap\fg^{\sigma'}$, which is a subalgebra
of type $G_2$.  Since the only intermediate subalgebra between $G_2$
and $\fg$ has type $B_3$, we see that $\fg^\rho$ is simple.  It
follows from Theorem~\ref{t:ssAuto} that $\rho$ is semisimple and inner;
hence, $\fg^\rho$ has rank equal to the rank of $\fg$.  Thus,
$\fg^\rho=\fg$, 
$\rho=1$, and $\tau=\sigma'$.  This implies that
$\fg^\sigma=\fg^{\sigma'}$, which is a contradiction. Hence,  $\overline{\tau}=\overline{\sigma}$ and by the same argument as above, $\tau=\sigma$.

\subsubsection{} Finally, suppose $\overline{\tau}$ and $\overline{\sigma}$ have order $3$. 
Then,  $\fg^\sigma=\fg^{\sigma^2}$, and
$\overline{\tau}$ equals either $\overline{\sigma}$ or
$\overline{\sigma^2}$. By the same argument as above, $\tau=\sigma$ or $\tau=\sigma^2$. In particular, $\tau$ is pinned.  
This concludes the proof. 
 \qed

\subsection{Connectedness criterion for cyclic subgroups}\label{ss:connected}

Given an element $g \in G$, let $\overline{\langle g \rangle}\subseteq G$ denote the Zariski closure of the cyclic subgroup generated by $g$. In this subsection, we provide a criterion for determining when $\overline{\langle g \rangle}$ is connected. 

 Let $g = g_s g_u$ be the Jordan decomposition of $g$. 
 \begin{lem} We have an isomorphism of algebraic groups $\overline{\langle g \rangle} \simeq \overline{\langle g_s \rangle} \times \overline{\langle g_u \rangle}$. 
 \end{lem} 
\begin{proof} As $g = g_s g_u \in \langle g_s \rangle \times \langle g_u \rangle$, we conclude 
  $\langle g \rangle \subseteq \langle g_s \rangle \times \langle g_u \rangle$. Thus, 
  \[ 
    \overline{\langle g \rangle} \subseteq \overline{\langle g_s \rangle \times \langle g_u \rangle} 
    \simeq \overline{\langle g_s \rangle} \times \overline{\langle g_u \rangle},
      \]
  where the last isomorphism follows from \cite[\S II.2.1]{Borel}.

  Conversely, since an algebraic group contains the semisimple and unipotent parts of its elements, we have $(g^m)_u = (g_u)^m \in \overline{\langle g \rangle}$ and $(g^n)_s = (g_s)^n \in \overline{\langle g \rangle}$ for all integers $m, n$. Hence, $(g_s)^n (g_u)^m \in \overline{\langle g \rangle}$ for all $m, n$. Thus, we can identify  $\langle g_s \rangle \times \langle g_u \rangle$ with a subgroup of  $\overline{\langle g \rangle}$.
  Taking closures, we obtain the reverse inclusion.
\end{proof}

\begin{exam} Let $G=\SL_2$ and  $g=\left(\begin{matrix} x & 1 \\ 0 & x^{-1}\end{matrix}\right)$ with $x\in \bC^\times$. If  $x$ is a primitive $e$-th root of unity, then $\overline{\langle g \rangle} \simeq \mu_e \times U$, where $U=[B,B]$ is the group of unipotent upper triangular matrices; otherwise, $\overline{\langle g \rangle} = B$. One way to prove this is to note that $\overline{\langle g \rangle}$ is either one or two dimensional and use the fact that the only connected one-dimensional linear algebraic groups over $\bC$ are $\mathbb{G}_m$ and $\mathbb{G}_a$. 
\end{exam} 

\subsubsection{} 
To state our connectedness criterion, we need to introduce some notation. 
Choose $t \in T$ in the conjugacy class of $g_s$, and define the \emph{group of character values} of $g$ by
\[
X^*(g) := \{ \lambda(t) \mid \lambda \in \Hom(T, \mathbb{C}^\times) \} \subseteq \mathbb{C}^\times.
\]
Note that $X^*(g)$ is a finitely generated subgroup of $\mathbb{C}^\times$ which depends only on the conjugacy class of $g_s$. When $G \subseteq \GL_n(\mathbb{C})$, the group $X^*(g)$ is generated by the eigenvalues of $g$.

An analogous construction exists in the Lie algebra setting. Let $A \in \mathfrak{g}$, and choose $B\in \ft$ in the adjoint orbit of the semisimple part $A_s$. Define 
\[
X^*(A) := \{ d\lambda(B) \mid \lambda \in \Hom(T, \bC^\times) \} \subseteq \mathbb{C}.
\]
Here, $d\lambda: \ft \ra \bC$ denotes the differential of the character $\lambda$. Again,  $X^*(A)$ depends only on the adjoint orbit of $A_s$. 

\begin{defe} \label{d:torsionFree} We say $g\in G$ is torsion-free if $X^*(g)$ is a torsion-free abelian group. Likewise,  $A\in \fg$ is torsion-free if  $X^*(A) \cap (\mathbb{Q} \setminus \mathbb{Z}) = \varnothing$.
\end{defe} 

Note that if $g\in \GL_n$, then $X^*(g)$ is torsion free if and only if $g$ has no eigenvalue which is a nontrivial root of unity. Similarly, for $A\in \mathfrak{gl}_n$,   $X^*(A)$ is torsion-free if and only if  $A$ has no rational non-integral eigenvalues. 

\begin{lem} \label{l:connected}
  \begin{enumerate}
    \item The algebraic group $\overline{\langle g \rangle}\subseteq G$ is connected if and only if $g$ is torsion-free. 
    
    \item The algebraic group $\overline{\langle \exp(2\pi i A) \rangle}\subseteq G$ is connected if and only if
   $A$ is torsion-free. 
  \end{enumerate}
\end{lem}

\begin{proof} By the previous lemma $\overline{\langle g \rangle} = \overline{\langle g_s \rangle} \times \overline{\langle g_u \rangle}$. Moreover, $ \overline{\langle g_u \rangle}$ is connected; in fact, it is trivial or isomorphic to $\mathbb{G}_a$. Thus, we may assume without loss of generality that $g\in T$.

Suppose $X^*(g)$ contains torsion.  Then there exists $\lambda\in X^*(T)$ such that $\lambda(g)$ is a primitive $e$th root of unity for some integer $e > 1$. By continuity, $\lambda(\overline{\langle g \rangle})$ is a finite group of order $e$, so $\overline{\langle g \rangle}$ is not connected.

  Conversely, assume $\overline{\langle g \rangle}$ is not connected. Since it is a diagonalisable group, it is isomorphic to $S\times F$, where $S$ is a torus and $F$ is a nontrivial finite abelian group.  We then have an isomorphism $X^*(\overline{\langle g \rangle}) \simeq X^*(S) \times F^\vee$, where $F^\vee$ is the Pontryagin dual. If $\nu(g) = 1$ for all $\nu \in F^\vee$, then $g \in S$, contradicting the assumption that $\overline{\langle g \rangle}$ is not contained in $S$. Therefore, there exists $\nu \in F^\vee$ such that $\nu(g)$ is a nontrivial root of unity. Since $\nu$ is the restriction of a character of $T$, it follows that $X^*(g)$ contains torsion.

 Part (2) follows  from (1), since the condition on $X^*(A)$ is equivalent to $X^*(\exp(2\pi i A))$ containing no roots of unity.
\end{proof}

\subsection{Irreducible and reductive subgroups} \label{ss:irreducible} The goal of this subsection is to show that certain subgroups of $G$ containing a lift of a Coxeter element are irreducible and reductive. 

Following Serre, we say that a closed subgroup $H \subseteq G$ is \emph{reducible} (in $G$) if there exists a proper parabolic subgroup $P \subseteq G$ such that $H \subseteq P$. Otherwise, $H$ is \emph{irreducible}. A subgroup $H \subseteq G$ is \emph{completely reducible} if whenever $H$ is contained in a parabolic subgroup $P \subseteq G$, there exists a Levi subgroup $L \subseteq P$ such that $H \subseteq L$. Note that every irreducible subgroup is, in particular, completely reducible. 

\begin{lem} \label{l:overGroup} 
Suppose $H\subseteq K\subseteq G$. If $H$ is irreducible in $G$, then so is $K$. 
\end{lem} 
\begin{proof} Indeed, if $K$ is reducible, then so is every subgroup of $K$. 
\end{proof} 

\subsubsection{} The following theorem, which goes back to  works of Serre \cite{Serre} and Borel--Tits \cite{BorelTits},  plays a foundational role in this subject. We refer the reader to \cite{Litterick} for further details.

\begin{thm} \label{t:BT}
A connected subgroup $H\subseteq G$ is completely reducible if and only if it is reductive. 
\end{thm}

\subsubsection{} \label{sss:NGT} We now show that $N_G(T)$ is an irreducible subgroup of
  $G$. To this end, we need a lemma. 
\begin{lem}  Let $P$ be a parabolic subgroup of $G$ with the
  unipotent radical $U$. Let $L$ be the unique Levi subgroup of $P$
  containing $T$. Then $N_G(T)\cap P =N_L(T)$. 
\end{lem}

\begin{proof}  Let  $n\in N_G(T)\cap P$, and write $n=l u$ for $l\in L$
  and $u\in U$.  Given any $t\in T$, we have
  \begin{equation*}  ntn^{-1}=l u t u^{-1} l^{-1} =\left( l t l^{-1}\right)\left( l
   \left[ (t^{-1}  u t) u^{-1}\right] l^{-1}\right),
  \end{equation*}
  with the first factor on the right in $L$ and the second in $U$.
  Since $ntn^{-1}\in T\subset L$, the factor in $U$ is the
  identity.  It follows that $t^{-1}  u t u^{-1} =1$, so $t^{-1} u t =
  u$ for all $t\in T$.  The only unipotent element in the centraliser
  of $T$ is the identity, so $u=1$ and $n\in L$.
\end{proof} 

\begin{cor} The subgroup $N_G(T)$ is irreducible in $G$.  
\end{cor} 

\begin{proof} Suppose $N_G(T)$ is contained in a proper parabolic $P=LU$ with $L$ containing $T$. The previous lemma states 
 that $N_G(T)=N_L(T)$. Thus, 
the  Weyl group $W_L$ of $L$ coincides with $W$.  This is a contradiction
 as $W_L$ is a proper parabolic subgroup of $W$.
\end{proof}

\subsubsection{Subgroups containing a lift of a Coxeter element} \label{sss:Coxeter}
Now, let $\theta\in N_G(T)$ be a finite-order lift of a Coxeter element $c\in W$; such a lift exists by a theorem of Tits \cite{Tits}.  Let $S\subseteq T$ be the smallest torus whose Lie algebra contains a fixed regular eigenvector of $c$. As noted in \cite[\S 13]{FG}, $\dim(S)=\phi(h)$, where $\phi$ is the Euler totient function. Let 
\[
H:=S \rtimes \langle \theta\rangle\subset G.
\]
Since $\theta$ has finite order,  $H$ is a closed subgroup of $G$, and $H^\circ=S$.

\begin{lem} \label{l:irreducible}
The subgroup $H$ is irreducible in $G$. 
\end{lem} 

\begin{proof}  Suppose $H$ lies in some  parabolic
  subgroup $P\subseteq G$. We first show that $P$ contains $T$.  Since
  $S$ contains a regular element, $T$ is the unique maximal
  torus in $G$ containing $S$.  However, $S$ is contained in a maximal
  torus in $P$, which is also maximal in $G$.  We conclude that
  $T\subseteq P$.

Now, let $L\subseteq P$ be the Levi
subgroup containing $T$. By the above lemma, $\theta\in N_L(T)$.  (Note this already shows that
$H\subseteq L$, so $H$ is completely reducible.)  Since $\theta$ has a primitive $h$th root of unity as an eigenvalue, 
 it follows from \cite[Theorem~3.4]{Springer} that $h$ divides a fundamental degree of $L$. By inspection (or
 by Theorem \ref{t:large} below), we must have $L=G$.
  Thus, $P=G$, and  $H$
 is irreducible.
\end{proof}

\begin{prop} \label{p:reductive} Let $K\subseteq G$ be a connected subgroup containing $H=S \rtimes \langle \theta\rangle$. Then $K$ is irreducible and reductive. 
\end{prop} 

\begin{proof} Since $H$ is irreducible, so is $K$. Theorem \ref{t:BT} then implies that $K$ is reductive. 
\end{proof}

\subsection{Coxeter grading} \label{ss:CoxeterGrading}
In this subsection, we recall basic facts about the Coxeter grading and exponents.

  \subsubsection{} 
The \emph{Coxeter grading} of $\fg$ (associated to $(B,T)$) is the $h$-periodic grading
$\fg=\bigoplus_{i\in \Z/h\Z} \fg_i$ 
given by
\[
\fg_i=\begin{cases} \bigoplus_{
   \mathrm{ht}(\alpha)\equiv i} \fg_\alpha, &\text{for $i\not\equiv 0$},\\
\ft,&\text{for $i\equiv 0$}.
\end{cases}
\]

\subsubsection{} One can show that there exists $n \in N_G(T)$
such that \[
\fg_i= \{ v \in \fg \mid n\cdot v = \zeta_h^i v \},
\]
where $\zeta_h=e^{2\pi i/h}$. In this case, $n$ must be a lift of a Coxeter element.

\subsubsection{Exponents}
An integer $e \in \{1,2,\dots,h\}$ is called an \emph{exponent} of $\fg$ if
$\exp(2\pi i e/h)$
is an eigenvalue of a Coxeter element acting on $\ft$.
Equivalently, $e$ is an exponent if and only if $\fg_e \neq 0$ in the Coxeter grading.

\subsubsection{Centraliser of $N+E$} Let $\ft'\subseteq \fg$ denote the centraliser of the regular semisimple element $N+E$ defined in \S \ref{ss:Notation}. 
Then $\ft'$ is a Cartan subalgebra, and the restriction of the Coxeter grading to $\ft'$ induces a grading 
\[
\ft' = \bigoplus_{i \in \mathbb Z/h\mathbb Z} \ft'_i .
\]
The graded piece $\ft'_i$ is nonzero if and only if $i$ is congruent modulo $h$
to an exponent of $\fg$.
In particular, $\ft'_0 = 0$.
Moreover,  the space $\ft'_i$ is one-dimensional (when non-zero),
except in type $D_{2n}$, where $\ft'_{2n-1}$ has dimension $2$.

\subsubsection{} \label{sss:XY}
Each element of $\ft'_i$ has nonzero components only in root spaces of height
congruent to $i$ modulo $h$.
In particular, if $s \in \{1,\dots,h-1\}$, then every element of $\ft'_{-s}$ admits a unique
decomposition
$X+Y$, 
where $X$ is a sum of elements in root spaces of height $-s$, and
$Y$ is a sum of elements in root spaces of height $h-s$.

\subsubsection{} It is known (see \cite{RLYG}) that if $d$ is the period of a finite positive rank grading\footnote{A $d$-periodic grading has positive rank if $\fg_r$ contains a non-nilpotent element for some integer $r$ with $\gcd(r,d)=1$.} of $\fg$, then $d \leq h$, with equality $d = h$ if and only if the grading is conjugate to the Coxeter grading. This fact will be used in the proof of the following lemma.

\begin{lem} \label{l:Coxeter}
Let $b$ be a positive multiple of $h$.
Let $\theta \in G$ satisfy $\theta^b = 1$, and suppose there exists a nonzero element
$x \in \ft$ such that
\[
\theta \cdot x = \zeta_h x.
\]
Then $x$ is regular, $\theta$ belongs to $N_G(T)$, and the image of $\theta$ in
$W = N_G(T)/T$ is a Coxeter element.
\end{lem}

\begin{proof}
The element $\theta$ induces a  $\mathbb Z/b\mathbb Z$-grading on $\fg$:
\[
\fg = \bigoplus_{i \in \mathbb Z/b\mathbb Z} \fg_i,
\qquad
\fg_i = \{ v \in \fg \mid \theta \cdot v = \zeta_b^i v \},
\]
where $\zeta_b = e^{2\pi i/b}$.
Since $b$ is divisible by $h$, we may collapse this grading to obtain an
$h$-periodic grading
\begin{equation} \label{eq:grading}
\fg = \bigoplus_{j \in \mathbb Z/h\mathbb Z} \fg'_j,
\qquad
\fg'_j = \bigoplus_{i \equiv j \!\!\! \pmod h} \fg_i .
\end{equation}

By assumption, $x$ lies in $\fg'_r$ for some $r$ coprime to $h$.
In particular, the grading \eqref{eq:grading} has positive rank.
By the above discussion, it is conjugate to the Coxeter grading. It follows from the properties of the Coxeter grading that every element of $\fg'_r$ is
 either regular semisimple element or nilpotent. 
As $x\in \ft$, it follows that it is regular.

Now, let $S$ be the smallest subtorus of $T$ whose Lie algebra contains $x$.
Since $\theta \cdot x = \zeta_h x$, $\theta$ normalises $S$, and therefore
normalises its centraliser $C_G(S)=T$.
Thus, $\theta$ belongs to $N_G(T)$ and it is a lift of a Coxeter element. 
\end{proof}

\subsection{Loop algebras and loop groups} \label{ss:loop} Let $\cK:=\bC(\!(t)\!)$ denote the field of Laurent series, $\fg(\cK)=\fg\otimes \cK$ the loop algebra, and $G(\cK)$ the loop group.

\subsubsection{Coxeter torus} \label{sss:CoxTorus}
Let 
\[
\omega_{-1} := N + t^{-1}E \in \mathfrak{g}(\mathcal{K}).
\]
This element is regular semisimple. Its centraliser $\mathcal{C}$ is a maximal anisotropic torus in $G(\cK)$. Under the bijection 
\[
\{\textrm{conjugacy classes of maximal tori in $G(\cK)$}\} \longleftrightarrow \{\textrm{conjugacy classes in $W$}\}
\]
 $\mathcal{C}$ corresponds to the Coxeter class. We refer to $\mathcal{C}$ as the \emph{Coxeter torus} (associated with our choice of affine pinning), and to its Lie algebra $\mathfrak{c}$ as the \emph{Coxeter Cartan subalgebra}.

\subsubsection{Iwahori grading} The ``standard Iwahori grading" on   $\fg[t,t^{-1}]$ is the $\frac{1}{h}\Z$ grading defined by  
\begin{equation}\label{eq:Iwahori}
\fg_I(m):=t^m \ft,\quad  \fg_I(m-s/h):=\bigg(\bigoplus_{\mathrm{ht}(\alpha)=h-s} t^{m-1}\fg_\alpha\bigg) \;\oplus\; \bigg(\bigoplus_{\mathrm{ht}(\alpha)=-s} t^m\fg_\alpha\bigg),\quad  m\in \mathbb{Z}, \quad s\in \{1,\dots, h-1\}.  
\end{equation} 
Note that $\omega_{-1}\in \fg_I(-1/h)$.

\subsubsection{Restriction to the Coxeter Cartan subalgebra} 
Setting  
\[
\fc(m-s/h):=\fc\cap \fg_I(m-s/h),
\]
we obtain a grading on  $\fc\cap
\fg[t,t^{-1}]$. Using the notation introduced in \S \ref{sss:XY} for the elements
of $\ft'$, the map
\[
\ft'_{-s}\to\fc(m-s/h),\qquad X+Y\mapsto t^{m}(X+t^{-1}Y). 
\]
is a linear isomorphism.

\subsubsection{Action of pinned automorphisms}\label{sss:action} The group
$\Aut(G, \cP)$ preserves the height of roots; in particular, it fixes
$N$ and acts on $\fg_I(j/h)$ for all $j\in \Z$. If we assume that $G$
is not of type $A_{2n}$, then $\Aut(G, \cP)$ also fixes $E$, and hence
the element $N+t^{-1}E$. It follows that, under this hypothesis,
$\Aut(G, \cP)$ preserves $\fc$ and acts on each graded piece
$\fc(j/h)$. We now describe this action explicitly in each type.

As the grading is $h$-periodic, we may assume
$j \in \{0,1,\dots,h-1\}$. First, consider the case of a pinned
involution $\sigma$. If $\fc(j/h)$ is one-dimensional,
then $\sigma$ acts trivially on it when $j$ is an exponent of
$\fg^{\sigma}$ and acts by $-1$ otherwise. In type $D_{2n}$ with
$n \geq 2$, and for $j = 2n-1$, the two-dimensional space $\fc((2n-1)/(4n-2))=\fc(1/2)$
decomposes as the direct sum of the two distinct irreducible
representations of $\langle \sigma \rangle \simeq \mathbb{Z}/2$.

It remains to consider the case of an automorphism $\tau$ of order
$3$ in type $D_4$. In this situation, $\tau$ acts trivially when
$j=1$ or $j=5$, while for $j=3$, the space $\fc(3/6)=\fc(1/2)$ decomposes as
the direct sum of the two nontrivial irreducible representations of
$\langle \tau \rangle$.

The above discussion shows that in all types (except $A_{2n}$), the
action of $\Aut(G, \cP)$ on $\fc(j/h)$ is trivial when $(j,h)=1$.

\subsubsection{} \label{s0}
Define $s_0$ to be $2$, $n-1$, $4$, and $3$ in types $A_{2n-1}$, $D_n$,
$E_6$, and $B_3$, respectively.  Except in type $B_3$, $s_0$ is the smallest
positive  integer $k$ such that $\Aut(G, \cP)$ does not act trivially on $\fc(-k/h)$. In type $B_3$, a similar statement holds if we consider automorphisms of a group $G'\supset G$ of type $D_4$. 

Note that Theorem \ref{t:large} below implies that $s_0$ is the smallest exponent 
that either does not appear, or---in type
$D_{2n}$---appears with smaller multiplicity, in any proper reductive
subgroup of maximal degree in $G$.

\subsection{Extended affine Weyl group} \label{ss:Weyl} In this subsection, we review key properties of the affine Weyl group and its variants, focusing on their actions on $\ft$.

\subsubsection{} 
Recall that we have inclusions 
\[
Q^\vee \subseteq X_*(T)\subseteq P^\vee,
\]
where $Q^\vee$, $X_*(T)$, and $P^\vee$ denote the coroot lattice, cocharacter lattice, and coweight lattice, respectively. Accordingly, we have three variants of the affine Weyl group, namely, the usual affine Weyl group, the extended affine Weyl group, and the extended affine Weyl group for the adjoint group: 
\[
W_\mathrm{aff}:=W\ltimes Q^\vee \quad  \subseteq \quad \widetilde{W}:=W\ltimes X_*(T)\quad  \subseteq\quad  \widetilde{W}_{\mathrm{ad}}:=W\ltimes
P^\vee.
\]

\subsubsection{} The relevance of these groups in our story emerges from the following lemma. 
Let 
\[\exp(2\pi i -) \colon \ft \to T
\]
 denote the exponential map.

\begin{lem}\label{l:WeylConjugate}
Let $S,S' \in \ft$. Then $S$ is $\widetilde{W}$-conjugate to $S'$ if and only if
$\exp(2\pi i S)$ is $W$-conjugate to $\exp(2\pi i S')$.
\end{lem}

\begin{proof}
Recall that the kernel of the exponential map $\exp(2\pi i -)\colon \ft \to T$
can be identified with $X_*(T) \subset \ft$.
By definition, $\exp(2\pi i S)$ is $W$-conjugate to $\exp(2\pi i S')$
if and only if there exists $w \in W$ such that
\[
\exp\bigl(2\pi i (w\cdot S - S')\bigr) = 1 .
\]
This holds if and only if $w\cdot S - S' \in X_*(T)$. On the other hand, the extended affine Weyl group $\widetilde{W} = X_*(T) \rtimes W$
acts on $\ft$ by
\[
(\lambda,w)\cdot S = w\cdot S + \lambda , \qquad \lambda \in X_*(T),\ w \in W .
\]
Thus the condition $w\cdot S - S' \in X_*(T)$ is equivalent to $S$ and $S'$ lying in the same
$\widetilde{W}$-orbit. \end{proof}

\subsubsection{Fundamental domains} Let  $\ft_\bR=X_*(T)\otimes_\bZ\bR$ be  the standard apartment and $\cA\subset \ft_\bR$  the fundamental alcove. Then $\overline{\cA}$ is a fundamental domain for the action of $W_{\mathrm{aff}}$ on $\ft_\bR$. 

Next, let $\Delta'\subset\Delta$ be the set of minuscule simple roots, i.e., those simple roots which appear with multiplicity one in the decomposition of the highest root $\alpha_0$. Then,  the closure of 
\begin{equation*} \cA'=\{\mu \in\cA\mid \langle \mu,
  \alpha+\alpha_0\rangle< 1 \text{ for all $\alpha\in\Delta'$}\} \subset \ft_{\bR}
  \end{equation*}
    is a fundamental domain for
  $\widetilde{W}_{\mathrm{ad}}$ \cite{KomPre}. For a picture of $\cA'$ in types $A_2$ and $B_2$, see \cite{SecoGarnierNeeb25}.

 \subsubsection{The set $\cF$} \label{sss:SetF}
Let $\cF:=\cA'+i\ft_\bR$. This is an open subset of $\ft$ in the analytic topology.

  \begin{lem}\label{l:SetF} The set $\cF$  is preserved by the action  $\Aut(G, \cP)$ and contains at
    most one element in any $\widetilde{W}$-orbit.  
  \end{lem}
  \begin{proof} 
  Indeed, $\cA$ is defined in terms of inequalities
  involving the simple roots and the highest root.  Since every
  $\sigma\in \Aut(G, \cP)$ permutes the simple roots and fixes the
  highest root, it preserves $\cA$.  Moreover, pinned
  automorphisms permute the minuscule simple roots, so $\sigma$ also
  preserves $\cA'$.

   Next, suppose
    \[
    \mu'+i\nu'=\widetilde{w}\cdot (\mu+i\nu),\qquad \mu,\mu'\in\cA',\quad \nu,\nu'\in\ft_\bR,\quad \widetilde{w}\in \widetilde{W}.
    \]
      We then have
    $\widetilde{w}\cdot \mu=\mu'$. 
    Since $\overline{\cA'} \subset \ft_{\bR}$ is a fundamental domain for $\widetilde{W}_{\mathrm{ad}}\supseteq \widetilde{W}$, $\cA'$ intersects each $\widetilde{W}$-orbit in at most one point. 
This implies $\mu=\mu'$ and $\widetilde{w}$ is the identity.  It
    follows that $\nu=\nu'$.
  \end{proof}

\begin{cor}\label{c:SetF} Let $\sigma \in \Aut(G, \cP)-\{1\}$ and  $S\in \cF-\cF^\sigma$. Then $S$ and $\sigma(S)$ are not $\widetilde{W}$-conjugate. 
\end{cor} 

\begin{proof} 
Since $\sigma(S)$ and $S$ are distinct elements of $\cF$, the previous
lemma implies that they are not in the same $\widetilde{W}$-orbit.  
\end{proof}

\subsubsection{More on the complement $\cF - \cF^\sigma$}
Let $\ft^\sigma \subset \ft$ denote the fixed-point subspace of $\sigma$.
Since $\sigma \neq 1$, the subspace $\ft^\sigma$ is a proper closed linear
subspace of $\ft$. Consequently,
\[
\cF^\sigma := \cF \cap \ft^\sigma
\]
is a proper closed subset of $\cF$, and its complement
$\cF - \cF^\sigma$ is a nonempty open subset of $\cF$ in the analytic
topology.

\subsubsection{} 
As an example, consider type $A_{n-1}$.  Here,
\[\cA=\left\{\diag(a_1,\dots,a_n)\mid \text{$a_i\in\bR$ for all $i$, 
    $a_1>\dots>a_n$, $a_1-a_n<1$,  and\, 
    $\sum_{i=1}^n a_i=0$}\right\}.\]
  Since all simple roots are minuscule, $\cA'$ is the subset of $\cA$
  satisfying the additional inequalities
  \[ a_1+a_i-a_{i+1}-a_n<1\qquad \text{for $1\le i<n$}.\]
   Finally, the pinned involution maps $\diag(a_1,\dots,a_n)$ to $-\diag(a_n,\dots,a_1)$.


\section{Reductive subgroups of maximal degree} \label{s:large} 
In this section, we investigate reductive subgroups of maximal degree in simple groups. We begin in \S\ref{ss:exponents} by examining a result about fundamental degrees of a reductive subgroup dividing those of the ambient group. This  naturally leads to \S\ref{ss:large}, where we state and prove the main theorem of the section: a classification of reductive subgroups of maximal degree. Finally, in \S\ref{ss:fixed}, we present a refinement of our results needed for the case of $D_4$ and prove Proposition \ref{p:PinnBig}.

\subsection{Fundamental degrees of reductive subgroups} \label{ss:exponents}
 The following result is likely known but we were unable to find a reference in the literature.

\begin{prop}\label{p:degrees} 
Let $K$ be a (connected, closed) reductive subgroup of $G$. Let $b$ be a fundamental degree of $K$. Then $b$ divides a fundamental degree of $G$. In particular, $b \leq h(G)$.
\end{prop}

 As an immediate corollary, we obtain another characterisation of reductive subgroups of maximal degree.

\begin{cor} Let $K$ be a reductive subgroup of $G$, and suppose one of the fundamental degrees of $K$ equals $h(G)$. Then $K\subseteq G$ is a reductive subgroup of maximal degree, i.e., $h(K)=h(G)$. 
\end{cor}

For a uniform proof of the above proposition using formal connections, see \S \ref{ss:degree}.  One can also prove the proposition by a case-by-case analysis. Indeed, we may assume that $K$ is semisimple, and it suffices to consider the case where $K$ is a \emph{maximal} semisimple subgroup of $G$. Such subgroups were classified by Dynkin \cite{Dynkin}, and one can verify the proposition by consulting the tables in his classification. It would be interesting  to find a uniform group-theoretic proof. 

\subsubsection{} We can give a uniform group-theoretic proof of Proposition \ref{p:degrees} when $K\subseteq G$ has maximal rank, i.e., contains a maximal torus of $G$. Without loss of generality, we may assume $T\subseteq K$. Let $\fk$ denote the Lie algebra of $K$ and $W_K$ its Weyl group.  Since $b$ is a fundamental degree of $K$,  \cite[Theorem~3.4]{Springer} implies that there exists $A \in\ft $ and $w \in W_K$ such that
\[
w \cdot A = e^{2\pi i/b} A.
\]

Now,  $W_K=N_K(T)/T$ is a subgroup of $W_G=N_G(T)/T$. Thus, the element $w$ also lies in $W_G$. Appealing again to \cite[Theorem~3.4]{Springer}, we conclude that $b$ divides one of the fundamental degrees of $G$. 
\qed

\subsection{The classification theorem} \label{ss:large}
The following is a more precise version of Theorem \ref{t:largeIntro}. 

\begin{thmbis}{t:largeIntro}\label{t:large} Let $G$ be a simple group. In all types except $A_{2n}$ and $B_3$, the following are equivalent: 
\begin{enumerate} 
\item $K\subseteq G$ is a reductive subgroup of maximal degree. 
\item $K$ is $G$-conjugate to $G^{\sigma, \circ}$ for some $\sigma\in \Aut(G, \cP)$. 
\end{enumerate} 
 In type $A_{2n}$,  there are no proper reductive subgroups of maximal degree.\\ 
In type $B_3$, there is a proper reductive subgroup $K$ of maximal degree of type $G_2$. This subgroup is conjugate to $(G')^{\tau, \circ}$ where $G'$ is a group of type $D_4$ containing $G$ and $\tau\in \Aut(G', \cP')$ is an element of order $3$.  The pinning  $\cP'$ is chosen so that $(G')^\sigma=G$ for an involution $\sigma\in \Aut(G', \cP')$

A reductive subgroup of maximal degree is simple, and aside from the case $B_3\hookrightarrow D_4$,  unique up to $G$-conjugacy. Moreover, if
$K$ is proper, then $\fk\hookrightarrow \fg$ is maximal, except in the case $G_2\hookrightarrow D_4$ in which case it is maximal in a subalgebra of type
$B_3$. The root systems for all possible pairs $(G,K)$ are given in
Table \eqref{ta:main}.

\begin{table}[h]
\renewcommand{\arraystretch}{1.2} 
\centering
\begin{tabular}{|l|l|} 
\hline
$G$  &  $K$ \\ \hline
$A_1$ & $A_1$ \\ \hline
$A_{2n-1},\ n \geq 2$ & $A_{2n-1}$, $C_n$ \\ \hline
$A_{2n},\ n \geq 1$ & $A_{2n}$ \\ \hline
$B_3$ & $B_3,\ G_2$ \\ \hline
$B_n,\ n \geq 4$ & $B_n$ \\ \hline
$C_n,\ n \geq 2$ & $C_n$ \\ \hline
$D_4$ & $D_4,\ B_3,\ G_2$ \\ \hline
$D_n,\ n \geq 5$ & $D_n,\ B_{n-1}$ \\ \hline
$E_6$ & $E_6,\ F_4$ \\ \hline
$E_7$ & $E_7$ \\ \hline
$E_8$ & $E_8$ \\ \hline  
$G_2$ & $G_2$ \\ \hline
$F_4$ & $F_4$ \\ \hline
\end{tabular}
\vspace{0.8em} \caption{Root systems of reductive subgroups of maximal degree}
\label{ta:main} 
\end{table}
\end{thmbis}

The remainder of this subsection is devoted to the proof of the theorem, which proceeds via a case-by-case analysis. (It would be desirable to have a type-independent proof.)  

We begin by considering the statement at the level of Lie algebras. By assumption, $[\fk,\fk]$ contains a simple subalgebra $\fk'$ with Coxeter number $h = h(G)$. We first classify the possible types of $\fk'$ in each case. Then, we prove that $\fk'=\fk$ and discuss the implications at the group level. We shall use the  inclusions $\fg^\sigma \hookrightarrow \fg$ classified in Table \ref{ta:pinned}. 

 \subsubsection{The case $\fg=\Sl_n$} By assumption, $h(\fk')=n$; thus, $\fk'$ can have the following types: 
\begin{itemize}
\item $B_{n/2}$ or $D_{n/2+1}$, $n$ even. This is impossible, because Lie algebras of these types have no faithful $n$-dimensional representations. 
\item $C_{n/2}$, $n$ even. This is the embedding $\Sp_n \hookrightarrow \Sl_n$. 
\item Exceptional types. In this case, the exceptional subalgebra would have  Coxeter number $n$ and a faithful $n$-dimensional representation. Inspecting Coxeter numbers and dimensions of representations of Lie algebras of exceptional types shows that this is impossible. 
\end{itemize} 
 
 \subsubsection{The case $\fg=\So_{2n+1}$ with $n\geq 2$} By assumption, $h(\fk')=2n$; thus,  $\fk'$ can have the following types: 
\begin{itemize}
\item $A_{2n-1}$. This  is impossible because 
\[
\dim(\fk') = \dim(\Sl_{2n})= 4n^2-1> 2n^2+n=\dim(\fg). 
\]
\item $C_n$. In this case, $\dim(\fk')=\dim(\fg)$ so if $n>2$, this is impossible. (Note that for $n=2$, $\Sp_{4}=\So_5$.) 
\item $D_{n+1}$. This is also impossible because $\So_{2n+2}$ has no faithful representation of dimension less than $2n+2$.
\item Exceptional types. 
In this case, the exceptional subalgebra would have Coxeter number $2n$ and a faithful $2n+1$-dimensional representation. By inspection, we see that the only possibility is the embedding $\fg_2\hookrightarrow \So_7$.
\end{itemize}

 \subsubsection{The case $\fg=\Sp_{2n}$} In this case, $\fk'$ is also of maximal degree in $\Sl_{2n}$. Thus, our results for type $A$ imply that the only possibility is $\fk=\fg$.

\subsubsection{The case  $\fg=\So_{2n}$, $n\geq 4$} By assumption, $h(\fk')=2n-2$; thus, $\fk'$ can have the following types:
\begin{itemize}
\item $A_{2n-3}$. A dimension count shows that this is impossible.
\item $B_{n-1}$. This corresponds to the embedding $\So_{2n-1}\hookrightarrow \So_{2n}$. 
\item $C_{n-1}$. This is impossible because an embedding
  $\Sp_{2n-2}\hookrightarrow \So_{2n}$ does not exist.  To see this,
  suppose we have such an embedding.  Then we would have a
  decomposition $\C^{2n}=V\oplus W$ of $\Sp_{2n-2}$-modules, where $V$
  is the standard representation and $W$ is a two-dimensional trivial
  representation.  (This is the unique faithful representation of
  $\Sp_{2n-2}$ of dimension $2n$.)  If $\langle,\rangle$ is the
  underlying non-degenerate symmetric bilinear form, we obtain
  \[
  \langle Xv,w\rangle=-\langle v,Xw\rangle=-\langle v,0\rangle=0,\qquad \forall \, X\in \Sp_{2n-2}, \quad v\in V, \quad w\in W.
  \]
    By irreducibility of
  $V$, we conclude that $W=V^\perp$.  It follows that
  \[
  \Sp_{2n-2}\hookrightarrow \So_{2n-1}\oplus\So_2,
  \]
   which is
  impossible for dimension reasons.

\item Exceptional types. In this case, $\fk'$ would have Coxeter number $2n-2$ and a faithful $2n$-dimensional representation. By inspection, the only possibility is the embedding $\fg_2\hookrightarrow \So_8$. 

 \end{itemize}

\subsubsection{The case $\fg=E_6$}  In this case, the only possibilities for $\fk'$ are $A_{11}, B_6, C_6, D_7, F_4, E_6$. But dimension counts rule out the first four, and the fifth one corresponds to the embedding $F_4\hookrightarrow E_6$.

 \subsubsection{Remaining exceptional cases} Next, suppose $\fg$ is of type $E_7$ (resp. $E_8$). Then the only
 possibilities for $\fk'$ are  $A_{17}, B_9, C_9, D_{10}, E_7$
 (resp. $A_{29}, B_{15}, C_{15}, D_{16}, E_8$). But all the
 non-trivial possibilities are ruled out because their rank is greater
 than $7$. Finally, if $\fg$ is of type $G_2$ (resp. $F_4$), then
 $\fk'=\fg$ because of the above results for $B_3$ (resp. $E_6$). 
 
 This concludes the classification of possible $\fk'$ in each type. 

\subsubsection{Simplicity} 
We now prove that $\fk=\fk'$.  Without
 loss of generality, we can assume that $\ft\cap\fk$
 and $\ft\cap\fk'$ are Cartan subalgebras for $\fk$ and $\fk'$
 respectively.  If $\fk\ne\fk'$, then either $[\fk,\fk]=\fk'\oplus\mathfrak{m}$
 with $\mathfrak{m}$ semisimple or $\fk=\fk'\oplus\fz$, where $\fz$ is the
 nontrivial centre of $\fk$.  In either case, the centraliser of
 $\fk'$ in $\fk$ contains an element of $\ft-\fk'$.  However,
 in each case of $\fk'$ found above, it is easy to check that every nonzero element of $\ft$
 acts nontrivially on some root space of $\fk'$.  (The most
 complicated case is the inclusion of the Lie algebra of type $F_4$ in
 $E_6$, but it is well known that this subalgebra is maximal in
 $E_6$.)  This contradiction shows that $\fk$ must be simple. 
  
\subsubsection{Relationship to fixed points of pinned automorphisms}  
The above discussion shows that in each type except $B_3$ and $A_{2n}$,
\[
\textrm{$\fk\subseteq \fg$ is a reductive subalgebra of maximal degree} \quad \iff \quad \textrm{$\fk=\fg^{\sigma}$ for some $\sigma\in \Aut(G, \cP)$}.
\]
 If $\fg$ is of type $B_3$, then we have the additional case of $\fk = {\fg'}^{\tau}$, where  $\tau$ is a pinned automorphism of order $3$ of a Lie algebra $\fg'\supseteq \fg$ of type $D_4$. Finally, we have seen that there are no reductive subalgebras of maximal degree in type $A_{2n}$. 
By the discussion of \S \ref{sss:fixedPinn}, except for $G_2\hookrightarrow D_4$, reductive subalgebras $\fk\subseteq \fg$ of maximal degree are maximal. (We also prove this directly below.) In addition, except for the inclusion $B_3\hookrightarrow D_4$, the reductive subalgebra $\fk\subseteq \fg$ of maximal degree is unique up to conjugation.

\subsubsection{The group setting} Moving to the setting of algebraic groups, we claim that the reductive subgroup $K\subseteq G$ is  determined, up to conjugacy, by its root system. Indeed, let $H\subseteq G$ be a (connected) reductive subgroup of maximal degree with the same root system as $K$.  We already know that $H$ and $K$ are both simple,  and their Lie algebras are $G$-conjugate, i.e., there exists $g\in G$ such that $\Ad_g \fh=\fk$.  Thus, $gHg^{-1}$ and $K$ are two connected subgroups of $G$ with the same Lie subalgebra.  Since the map $\Lie(-)$ from connected algebraic subgroups of $G$ to Lie subalgebras of $\fg$ is injective, it follows that $gHg^{-1}=K$.

This concludes the proof of the theorem. \qed

\subsubsection{Irreducibility}  \label{sss:Irreducible}
Here, we show directly that subalgebras $\fk\subseteq \fg$ of maximal degree are irreducible. 
Suppose that $\fk\subseteq \fp$, where $\fp$ is a proper parabolic
subalgebra.  Since $\fk$ is reductive, it is completely
reducible, by Theorem \ref{t:BT}.  Thus, $\fk$ is a subalgebra of a proper Levi subalgebra
 $\fl\subseteq \fp$.  By Proposition \ref{p:degrees}, $h(\fk)\le h(\fl)$. On the other hand, 
as $\fk$ has maximal degree, it follows that $h(\fg)=h(\fk)\le h(\fl)$.  However, Table \ref{ta:main} implies  that $h(\fl)<h(\fg)$ for every proper Levi subalgebra.  This is a
contradiction. Thus, $\fk$ is irreducible. 

\subsubsection{Maximality} \label{sss:maximal}
We now give a direct proof that proper subalgebras $\fk\subseteq \fg$ of maximal degree are maximal, except for $G_2\hookrightarrow D_4$. Assume we are not in type $B_3$, so a proper reductive subalgebra of
maximum degree is of the form $\fg^\sigma$ for some nontrivial pinned
automorphism $\sigma$.  Let $\fk\subseteq\fg$ be a subalgebra containing
  $\fg^\sigma$.  Since
  $\fg^\sigma$ is irreducible,  so is $\fk$. By Theorem \ref{t:BT}, $\fk$ is reductive.
  Proposition \ref{p:degrees} implies $h(\fg^\sigma)\le h(\fk) \le
h(\fg)$. Since $h(\fg^\sigma)=
h(\fg)$, $\fk$ is reductive of maximal degree.  The statements about
maximality now follow from Table~\ref{ta:main}.  The case of $B_3$
follows from the result for type $D_4$.

\subsection{Fixed point subalgebras} \label{ss:fixed} In the proof of our main theorem (Theorem \ref{t:mainIntro}), we will need the following result about fixed point subalgebras.

\begin{prop}\label{p:pinned} 
Suppose $\fg^\tau$ contains a reductive subalgebra $\fk \subseteq \fg$
of maximal degree. Then $\tau$ is a pinned automorphism.
\end{prop}

\begin{proof}\footnote{We thank the contributors to the following  post:\\
\url{https://mathoverflow.net/questions/499355/fixed-point-subalgebras-of-automorphisms-of-d-4}}
By the discussion in \S \ref{sss:Irreducible}, $\fk$ is irreducible
  in $\fg$. Since $\fg^\tau$ contains $\fk$, it must also be
  irreducible in $\fg$. By Theorem \ref{t:BT}, $\fg^\tau$ is
  reductive. Next, Proposition \ref{p:degrees} implies
  \[
  h(\fk) \leq h(\fg^\tau)\leq h(\fg).
  \]
   Since $h(\fk)=h(\fg)$, these
  are all equalities. Thus, $\fg^\tau$ is a reductive subalgebra of
  maximal degree.

  If $\fg^\tau=\fg$, then $\tau=1$, and the result is immediate. Otherwise, $\fg^\tau$ is
  proper, and Theorem \ref{t:large} shows that $\fg^\tau$ is
  isomorphic to $\fg^\sigma$ for some pinned automorphism $\sigma$.
  By Proposition~\ref{p:PinnBig}, $\tau$ is a pinned automorphism.
\end{proof}


\section{Jordan form} \label{s:connections} In this section, we describe the Jordan canonical form (also known as the Levelt-Turrittin canonical form) of formal Coxeter connections and compute their differential Galois groups. In \S \ref{ss:connections} and \S \ref{ss:Auto} (resp. \S\ref{ss:formalConnections}) we review the essential definitions and constructions for connections (resp. formal connections). Using the notion of the slope of a formal connection, we prove Proposition~\ref{p:degrees} on fundamental degrees in \S\ref{ss:degree}. In \S\ref{ss:Jordan} and \S \ref{ss:localDiff}, we state and prove our main results concerning Jordan canonical forms and local differential Galois groups. 

In most of the first three subsections, we work in the general setting where $G$ is a connected complex linear algebraic group. In the last three subsections, we revert back to our usual conventions where $G$ denotes a simple group. 

\subsection{Connections on algebraic varieties}\label{ss:connections} Let $G$ be a connected complex linear algebraic group and $Y$ a smooth connected complex algebraic variety. A \emph{$G$-connection} on $Y$ is a pair $(\mathcal{E}, \nabla)$, where $\mathcal{E}$ is an algebraic (principal) $G$-bundle on $Y$, equipped with an algebraic connection $\nabla$. 
When the underlying bundle is clear from the context, we simply write $\nabla$ instead of $(\mathcal{E}, \nabla)$. If $\mathcal{E} = Y \times G$ is the trivial $G$-bundle, then $\nabla$ can be written as
\begin{equation*}
\nabla = d + A, \qquad A \in \Omega^1(Y, \mathfrak{g}),
\end{equation*}
where $d$ is the exterior derivative and $A$ is the connection one-form.
In general, after choosing a local trivialisation for $\mathcal{E}$, a $G$-connection admits a similar local description.

\subsubsection{Isomorphic connections} 
Two $G$-connections $(\mathcal{E}, \nabla)$ and $(\mathcal{E}', \nabla')$ are \emph{isomorphic} if there exists an isomorphism of $G$-bundles $\varphi: \mathcal{E} \to \mathcal{E}'$ such that
$\varphi^* \nabla' = \nabla$. 
If $\mathcal{E} = \mathcal{E}'$ is the trivial $G$-bundle, then every bundle automorphism $\varphi$ is given by an element $\mathtt{g} \in G(Y) = \mathrm{Hom}(Y, G)$.
Given connections $\nabla = d + A$ and $\nabla' = d + A'$ with $A, A' \in \Omega^1(Y, \mathfrak{g})$, $\varphi$ defines an isomorphism between $(\mathcal{E}, \nabla)$ and $(\mathcal{E}', \nabla')$ if and only if
\begin{equation*}
\nabla' = \mathtt{g}\cdot \nabla := d+ \mathtt{g} A \mathtt{g}^{-1} - \mathtt{g}^{-1} \, d \mathtt{g}.
\end{equation*}
This equivalence relation on $\Omega^1(Y, \mathfrak{g})$ is called \emph{gauge equivalence}.

\subsubsection{Tannakian perspective} 
Let $\pi_1^{\mathrm{diff}}(Y)$ denote the differential fundamental group of $Y$ \cite{Katz87}. To be precise, one should fix a base point $y \in Y$ and consider $\pi_1^{\mathrm{diff}}(Y, y)$, but we omit this from the notation for brevity.  

The category of $G$-connections on $Y$ is equivalent to the category of $G$-valued representations of $\pi_1^{\mathrm{diff}}(Y)$; in particular, two $G$-connections are  isomorphic if and only if their associated representations are $G$-conjugate. 

\subsubsection{Differential Galois group}
Given a $G$-connection $\nabla$ on $Y$ with corresponding representation $\rho_\nabla: \pi_1^{\mathrm{diff}}(Y) \to G$, the differential Galois group of $\nabla$ is defined by 
\[
G_\nabla:={\mathrm{Im}(\rho_\nabla)}.
\]
Note that $G_\nabla$ is automatically closed in $G$. Indeed, the representation $\rho_\nabla$ is algebraic and  factors through a finite type quotient to give a homomorphism $H\rightarrow G$ of algebraic groups, which has closed image. 

\subsection{Action of automorphisms on $G$-connections} \label{ss:Auto} 
We continue with the notation of the previous subsection. 
Let $\sigma\in \Aut(G)$. Then $\sigma$ induces an automorphism $\nabla \mapsto \sigma(\nabla)$ on the set of all $G$-connections on $Y$.
In local coordinates, if $\nabla = d + A$, then
\begin{equation*}
\sigma(\nabla) = d + \sigma(A).
\end{equation*}
The corresponding representation is given by  $\rho_{\sigma(\nabla)}=\sigma\circ \rho$.

\subsubsection{} \label{sss:Auto}
Recall that $\overline{\sigma}$ denotes the image of $\sigma$ under the canonical map $\Aut(G)\ra \mathrm{Out}(G)$. If $\sigma$ and $\tau$ are elements of $\Aut(G)$ with $\overline{\sigma}=\overline{\tau}$, then $\sigma(\nabla)\simeq \tau(\nabla)$. Indeed, in this case, 
 $\sigma=c_g\circ \tau$ for some $g\in G$. Thus, 
\[
\sigma(\nabla)\simeq g\cdot \tau(\nabla)\simeq \tau(\nabla). 
\]

\subsubsection{Justification of Remark \ref{r:forAll}} The previous paragraph together with the fact that  for $\overline{\tau}$ of order $3$, $\tau(\nabla)\simeq \nabla$ if and only if $\tau^2(\nabla)\simeq \nabla$, shows that one can replace every instance of ``for some" by ``for all" in Table \ref{ta:mainIntro}.

\subsubsection{Differential Galois group} Let $\nabla$ be a $G$-connection on $Y$.

\begin{lem}\label{l:auto}
\begin{enumerate} 
\item[(i)] $G_{\sigma(\nabla)}=\sigma(G_\nabla)$.
\item[(ii)] $\nabla \simeq \sigma(\nabla)$ if and only if $G_\nabla \subseteq G^{\tau}$ for some $\tau\in \Aut(G)$ satisfying $\overline{\tau}=\overline{\sigma}$. 
\end{enumerate} 
\end{lem} 

\begin{proof} 
For (i), note that  
\[
G_{\sigma(\nabla)} = {\operatorname{Im}(\rho_{\sigma(\nabla)})}
= {\operatorname{Im}(\sigma \circ \rho_\nabla)}
= {\sigma(\operatorname{Im}(\rho_\nabla))}
= \sigma(G_\nabla).
\]

For (ii), first suppose $\nabla \simeq \sigma(\nabla)$. Then the representations $\rho_\nabla$ and $\rho_{\sigma(\nabla)}$ are $G$-conjugate, i.e., there exists $g\in G$ such that 
\[
c_{g} \circ \rho_\nabla = \rho_{\sigma(\nabla)} = \sigma \circ \rho_\nabla.
\]
Thus, setting
$\tau := c_{g^{-1}} \circ \sigma$,  we obtain $\tau \circ \rho_\nabla = \rho_\nabla$. This implies $G_\nabla = {\mathrm{Im}(\rho_\nabla)}\subseteq G^\tau$.

Conversely, suppose $G_\nabla=\mathrm{Im}(\rho_\nabla) \subseteq G^\tau$ for some $\tau=c_g \circ \sigma\in \Aut(G)$. As noted in \S \ref{sss:Auto}, $\sigma(\nabla)$ and $\tau(\nabla)$ are then isomorphic. On the other hand, the representation $\rho_\nabla$ is fixed by $\tau$, i.e., $\tau \circ \rho_\nabla = \rho_\nabla$. 
Thus, $\nabla$ and $\tau(\nabla)$ are isomorphic. It follows that $\nabla$ and $\sigma(\nabla)$ are isomorphic. \end{proof}

\subsection{Formal connections} \label{ss:formalConnections}
For a positive integer $b$, let $\cK_b=\bC(\!(t^{1/b})\!)$, $\zeta_b=e^{2\pi i /b}$, and $\mu_b\subset \bC^\times$ be the group of $b$th roots of unity. We write $\cK$ for $\cK_1$. Let $\cD^\times:=\mathrm{Spec}\, \cK$ denote the formal punctured disk and $\cD^\times_b:=\mathrm{Spec}\, \cK_b$ its unique $b$-fold cover.

\subsubsection{} A \emph{formal $G$-connection} is a pair $(\cE, \nabla)$ consisting of a principal $G$-bundle on $\cD^\times$ and a connection $\nabla$ on $\cE$. Note that every principal $G$-bundle on $\cDt$ is trivial; after choosing a trivialisation, we can write $\nabla$ in the form $d + A \, dt$, where $ A \in \mathfrak{g}(\mathcal{K})$. Changing the trivialisation by $g\in G(\cK)$ transforms the connection matrix $A$ to $A'$ where 
\[
A':=g A g^{-1} - g^{-1} dg.
\]
For $ G = \mathrm{GL}_n $, the term $dg$ denotes the component-wise derivative (with respect to $t$) of the matrix $ g \in G(\cK)=G(\!(t)\!)$. For a general connected algebraic group $ G $, this expression can be defined by embedding $ G $ into some $ \mathrm{GL}_n $. The elements $A, A'\in \fg(\cK)$ are said to be  \emph{gauge equivalent} under the action of $G(\cK)$. 

\subsubsection{Regular versus irregular singularities} A formal connection is called \emph{regular singular} if it is gauge equivalent to $ d + A\,dt $, where $ A \in \mathfrak{g}(\mathcal{K}) $ has at most a simple pole. Otherwise, it is said to be \emph{irregular}.

\subsubsection{Local differential Galois group} As shown in \cite{Katz87}, we also have a local differential fundamental group $\pi_1^{\mathrm{diff}}(\cDt)$ such that the category $\mathrm{Conn}(\cD^\times)$ is equivalent to the category of representations of $\pi_1^{\mathrm{diff}}(\cD^\times)$. Thus, we can talk about ``the" representation $\rho_\nabla: \pi_1^{\mathrm{diff}}(\cDt)\ra G$ corresponding to a formal connection $\nabla$. 
As above, the image of this map is called the (local) \emph{differential Galois group} of $\nabla$ and denoted by $G_\nabla$.

\subsubsection{Slope} Suppose $G$ is reductive. Then a formal $G$-connection $ \nabla $ is, after a finite ramified extension $ \mathcal{K}_b / \mathcal{K} $, gauge equivalent to a connection of the form
\[
d+\Big(X t^{-\frac{a}{b}} + \text{higher order terms}\Big)\sdfrac{dt}{t},
\]
where $ X \in \mathfrak{g} $ is non-nilpotent and $a/b$ is a non-negative rational number, known as the \emph{slope} of $ \nabla $. For other characterisations of this notion, see \S \ref{ss:Strata} and \cite[\S 2]{CK}.

\subsection{Proof of Proposition \ref{p:degrees}} \label{ss:degree} 

Let $K$ be a connected reductive subgroup of $G$ with Lie algebra $\fk$.  Let $b$ be a fundamental degree of $K$.  Our goal is to show that $b$ divides a fundamental degree of $G$. 

\subsubsection{} 
Choose a maximal torus $S\subseteq K$, and let $W_K:=N_K(S)/S$ denote the Weyl group. By \cite[Theorem 3.4]{Springer}, there exists a non-zero element $A\in \fs:=\Lie(S)$ and $w\in W_K$ such that 
\begin{equation} \label{eq:Spring} 
w\cdot A = e^{-2\pi i/b} A.
\end{equation}

\subsubsection{} 
Now,  consider the $K$-connection on $\mathcal{D}_b^\times$ defined by 
\[
\nabla := d + (At^{-1/b}) \sdfrac{dt}{t}.
\]
By the descent criterion in \cite[\S 9]{BV}, 
 $\nabla$ is $K(\cK_b)$-gauge equivalent to a $K$-connection $\eta$ on $\cDt$ if and only if there exists $w\in W_K$ satisfying \eqref{eq:Spring}.\footnote{For a modern proof of this criterion, see \cite[\S 2]{JY}.}

\subsubsection{} 
Finally, observe that the slope of the $K$-connection $\eta$ is $1/b$. Let $\iota: K \to G$ denote the inclusion. Then the induced $G$-connection $\iota_*(\eta)$ on $\mathcal{D}^\times$ also has slope $1/b$. By \cite[§2]{CK}, $b$ divides a fundamental degree of $G$. This concludes the proof. 
\qed

\subsection{Jordan form}\label{ss:Jordan} 

The following is the main result of this section: 

\begin{thm}\label{thm:Jordan} 
Let $\nabla$ be a formal (Coxeter) $G$-connection of slope $r/h$, where $\gcd(r,h)=1$. Then there exists a positive integer $b$ divisible by $h$ and an element $\theta \in N_G(T)$ such that:
\begin{itemize}[label=$\bullet$]
    \item  $\theta^b = 1$, and the image of $\theta$ in the Weyl group $W$ is a Coxeter element $w$.
    \item $\nabla$ is $G(\mathcal{K}_b)$-gauge equivalent to a connection of the form
    $$
    d + (x_1 u^{r_1} + x_2 u^{r_2} + \cdots + x_n u^{r_n}) \frac{du}{u},
    $$
    where $u = t^{1/h}$, the integers $r_i$ satisfy $r_1 = -r < r_2 < \cdots < r_n < 0$, and each $x_i \in \mathfrak{t}$ is an eigenvector of $w$ with eigenvalue $e^{-2\pi i r_i/h}$.
\end{itemize}
\end{thm}

Before giving the proof, we recall the main tenets of the theory of Jordan canonical forms for formal connections; see, for example,
\cite[\S 9]{BV} or \cite{KW}. According to this theory, for every irregular singular formal $G$-connection
$\nabla$, there exists a positive integer $b$ and an element
$g \in G(\mathcal K_b)$ such that $g\cdot \nabla$ equals 
\begin{equation} \label{eq:Jordan}
 d + \bigl(x_1 v^{s_1} + \cdots + x_n v^{s_n} + y\bigr)\,\frac{dv}{v}.
\end{equation} 
Here, $v = t^{1/b}$, the integers $s_i$ satisfy
\[
s_1 < \cdots < s_n < 0,
\]
each $x_i \in \mathfrak t \setminus \{0\}$, and $y \in \fg$ commutes with every
$x_i$. Moreover, if we define
\[
\theta \;:=\; g(\zeta_b v)\, g(v)^{-1},
\]
then $\theta\in G$, and we have 
\begin{equation}\label{eq:eigenvector}
\theta^b = 1, \qquad
\theta \cdot x_i = \zeta_b^{s_i} x_i,
\quad \text{and} \quad
\theta \cdot y = y.
\end{equation}
The rational number $-s_1/b$ is  the slope of $\nabla$.

\begin{proof}[Proof of Theorem \ref{thm:Jordan}] Suppose the expression in \eqref{eq:Jordan} denotes the Jordan decomposition of $\nabla$. 
By assumption, the slope is $r/h$, i.e. $\lvert s_1/b \rvert = r/h$. Since $\gcd(r,h)=1$, we see that $h$ divides $b$. By Lemma \ref{l:Coxeter}, $x_1$ is regular, $\theta$ belongs to $N_G(T)$, and the image of $\theta$ in $N_G(T)/T=W$ is a Coxeter element, denoted by $w$.
Next, recall that $w$ has no eigenvalue equal to $1$ (see \cite[\S 3.16]{Humphreys}), hence $y = 0$. The eigenvalues of $w$ are $h$th roots of unity, so if we set $r_i := s_i h / b$, then each $r_i$ is an integer, and $x_i$ is an eigenvector of $w$ with eigenvalue $\zeta_h^{r_i}$. Letting $u := v^{b/h}$ yields the form stated in the theorem. 
 \end{proof}

For an alternate proof of this theorem using rational canonical forms, see \S \ref{ss:RatJor}.

\subsection{Differential Galois group} \label{ss:localDiff}
We continue using the notation of Theorem \ref{thm:Jordan}. Let $S$ be the smallest subtorus of $T$ whose Lie algebra contains $x_1, \dots, x_n$. 

\begin{cor} \label{c:localDiff} 
We have $G_\nabla = S \rtimes \langle \theta \rangle$. Moreover, $G_\nabla$ is irreducible in $G$. 
\end{cor} 

\begin{proof} It is well known (see \cite[Theorem 10.2]{SV}) that the differential Galois group of $\nabla$ fits into a split short exact sequence
\[
1 \to S \to G_\nabla \to \langle \gamma \rangle \to 1,
\]
where $S \subseteq G$ is the exponential torus and $\gamma \in G$ is the formal monodromy. The exponential torus is the smallest torus in $T$ whose Lie algebra contains $x_1, \dots, x_n$. Since the Jordan form of $\nabla$ has no residue term, the only contribution to the formal monodromy comes from the Galois action of $\operatorname{Gal}(\mathcal{K}_b / \mathcal{K})$, which corresponds to the action of $\langle \theta \rangle$ on $S$ determined by \eqref{eq:eigenvector}. Thus, $G_\nabla = S \rtimes \langle \theta \rangle$.

The last statement follows from Lemma \ref{l:irreducible}. 
\end{proof}


\section{Proof of the main theorem} \label{s:global} 

In this section, we prove Theorem \ref{t:mainIntro}. For ease of notation, let $K:=G_\nabla$. 

\subsection{$K$ is a reductive subgroup of maximal rank} 

We first show that $K$ is connected. Indeed, by a theorem of Gabber (\cite[\S 1.2.5]{Katz87}), the fact that $G_\nabla^{\mathrm{mono}}$ is connected implies that $K$ is connected. Note that Gabber worked with $\GL_n$-connections, but his theorem immediately generalises to $G$-connections.  One need only embed $G$ into $\GL_n(\C)$ for some $n$ and observe that $G_\nabla^{\mathrm{mono}}$ and $G_\nabla$ remain unchanged when viewing $\nabla$ as a $\GL_n$-connection.

Next, we prove that $K$ is reductive. Note that for every point $s\in S$, the differential Galois group of the formal connection $\nabla|_{\cDt_s}$ embeds in $K=G_\nabla$. If we restrict to a point $s$ with a Coxeter singularity,  Corollary \ref{c:localDiff} gives a description of the local differential Galois group $G_{\nabla_s}$. Proposition \ref{p:reductive} then implies that a connected group containing $G_{\nabla_s}$ is reductive. Thus, $K$ is reductive.    

Finally, we show that $K$ has maximal degree. 
As $G_\nabla=K$, we conclude that  $\nabla$ comes from a $K$-connection via the embedding $K \subseteq G$. Now  $\nabla|_{\cDt_x}$, thought of as a $K$-connection, also has slope $r/h$. By  \cite[\S 2]{CK},  $h$ divides a fundamental degree of $K$, and by Proposition \ref{p:degrees}, $h=h(K)$. Thus, $K$ is a connected reductive subgroup of $G$ of maximal degree.

\subsection{Type $A_{2n-1}$ for $n\ge 2$, $D_n$  for $n\geq 5$, or $E_6$} We now verify Table \ref{ta:mainIntro} in the cases where there is a unique nontrivial Dynkin automorphism. 

Suppose $\nabla\simeq \sigma(\nabla)$ for some non-inner automorphism $\sigma\in \Aut(G)$.  Lemma \ref{l:auto} then implies that $K\subseteq G^{\tau}$ for a (possibly different) non-inner automorphism $\tau \in \mathrm{Aut}(G)$. In particular, $K$ is a \emph{proper} subgroup of $G$. By Theorem \ref{t:large}, $K$ must be of type $C_n$, $B_{n-1}$, or $F_4$, respectively. 
  
 Conversely, suppose $\sigma(\nabla)\not\simeq \nabla$ for all non-inner automorphisms $\sigma\in \Aut(G)$. By Lemma \ref{l:auto}, $G_{\nabla}\not\subseteq G^{\tau}$ for every non-inner automorphism $\tau \in \Aut(G)$. In particular, $\not\subseteq G^{\sigma}$ for every nontrivial pinned automorphism $\sigma\in \PAut(G)$.  By Theorem \ref{t:large}, $G_\nabla=G$.

 \subsection{Type $D_4$} We now verify Table \ref{ta:mainIntro} for $D_4$.   Assume $\tau(\nabla)\simeq \nabla$
for some $\tau\in \Aut(G)$ with $|\overline{\tau}|=3$.  By Lemma \ref{l:auto},  $\fg_\nabla:=\Lie(G_\nabla)
\subseteq \fg^{\kappa}$ for some $\kappa \in \Aut(G)$ satisfying
$\overline{\kappa}=\overline{\tau}$. Proposition \ref{p:pinned} now
implies that $\kappa$ is pinned of order $3$. Thus, $\fg^\kappa$ is of type $G_2$. By Theorem \ref{t:large}, so is $G_\nabla$. 
 
Next, suppose $\#\{ \overline{\sigma} \, | \, \sigma(\nabla)\simeq \nabla,\,\, \sigma\in \Aut(G) \}=2$.  This is equivalent to the statement that 
$\sigma(\nabla)\simeq \nabla$ for a unique involution
$\sigma\in \Aut(G, \cP)$. By Lemma \ref{l:auto}, $\fg_\nabla\subseteq \fg^{\tau}$
for some $\tau\in \Aut(G)$ satisfying
$\overline{\tau}=\overline{\sigma}$. In particular, $\fg_\nabla$ is a
proper Lie subalgebra of $\fg$.  On the other hand, the uniqueness
assumption implies that $\kappa(\nabla)\not\simeq \nabla$ for any order
$3$ element $\kappa \in \Aut(G, \cP)$ because an order $2$ and an order
$3$ element of $S_3$ generate the whole group.  It follows that
$\rho(\nabla)\not\simeq \nabla$ for any pinned automorphism $\rho$
of order $3$.  By Lemma \ref{l:auto},
$\fg_\nabla \not\subseteq \fg^{\rho}$ for all such $\rho$.  As noted in \S \ref{sss:fixedPinn},
every subalgebra  $G_2\hookrightarrow \fg$ equals $\fg^{\rho}$ for such a
$\rho$. Thus, $G_\nabla$ does not have type $G_2$. In view of Theorem
\ref{t:large}, $G_\nabla$ has type $B_3$.

Finally, if neither of these cases occurs, then $\sigma(\nabla)\not\simeq \nabla$ for all nontrivial $\sigma\in \Aut(G,\cP)$. By Lemma \ref{l:auto}, $G_{\nabla}\not\subseteq G^{\tau}$ for every non-inner automorphism $\tau \in \Aut(G)$. By Theorem \ref{t:large}, $G_\nabla=G$.  
 
\subsection{Type $B_3$} Finally, we verify Table \ref{ta:mainIntro}
for type $B_3$. In view of Theorem \ref{t:large}, we need to consider
a group $G'$ of type $D_4$ such that $G=(G')^{\sigma, \circ}$ for some
element $\sigma\in \Aut(G)$ with $|\overline{\sigma}|=2$. We can then
think of $\nabla$ as a $G'$ connection and consider $G_\nabla$ as a subgroup of $G'$. 

The analysis above for $D_4$ now applies to show that the differential Galois group of $G_\nabla$ either has type $G_2$ or $B_3$, with the latter occurring precisely when 
$\tau(\nabla)\simeq \nabla$ for some (equivalently every) $\tau\in
\Aut(G')$ such that $\overline{\tau}$ has order
$3$. 

This concludes the proof of Theorem \ref{t:mainIntro}.

\qed


\section{Rational canonical form}  \label{s:rational}
In this section, we discuss the rational canonical form of formal Coxeter connections. We begin by stating our main theorem in \S \ref{ss:Rational}. To establish this result, we first recall key concepts from the theory of fundamental strata in \S \ref{ss:Strata}. The proof of our main theorem is then carried out in \S \ref{ss:Proof}. In \S \ref{ss:RatJor}, we show how the Jordan canonical form can be obtained from the rational canonical form. Finally, \S \ref{s:algo} provides an explicit algorithm for computing the rational canonical form of Coxeter connections.

\subsection{Statement of the main result} \label{ss:Rational} Recall that in \S \ref{ss:loop},  we defined the Coxeter torus $\cC$ in terms of the fixed affine pinning $(\cP, E)$. 
Let $r$ be a positive integer coprime to $h$.  

\begin{defe} 
The set of \emph{$\cC$-formal types  (or rational canonical forms) of depth $r/h$}, denoted by $\cA(\cC,r/h)$, is the subset
of $\bigoplus_{j=0}^r\fc(-j/h)$ with nonzero component in degree
$-r/h$. 

 The set of \emph{homogenous} formal types $\cA(\cC, r/h)_{\mathrm{hom}}$ is defined to be $\fc(-r/h)-\{0\}$. 
By the discussion of \S \ref{sss:action}, $\Aut(G, \cP)$ acts on the set $\cA(\cC, r/h)$ and fixes the subset $\cA(\cC, r/h)_{\mathrm{hom}}$.

\end{defe} 

The main result of this section is: 

\begin{thm}\label{t:CoxToral} Let $\nabla$ be a formal Coxeter $G$-connection of 
  slope $r/h$. Then $\nabla$ is $G(\cK)$-gauge equivalent to a
  connection of the form $d+A\dtt$ with  $A\in \cA(\cC,r/h)$.  Moreover, this
  representation is unique up to a free action of $\mu_h$ by scalar multiplication.
\end{thm}

We call $d+A\dtt$ a rational canonical form for $\nabla$. 
To prove this theorem, we need some recollections  from the  theory of fundamental strata and toral connections.

\subsubsection{}  In all types except $D_{2n}$, one can give a somewhat
  more explicit description of $\cA(\cC,r/h)$.  Recall the Cartan subalgebra $\ft'\subset \fg$ defined in \S \ref{ss:CoxeterGrading}. For each exponent $s$,
  choose a nonzero element $N_s+E_s\in\ft'_{-s}$, where $N_s$
  (resp. $E_s$) is a sum of root vectors of height $-s$ (resp. $h-s$).
  We will take $N_1=N$ and $E_1=E$.  Set
  \[
  \om_{-s}=N_s+t^{-1}E_s\in \fc(-s/h).
  \]
    Then,
\[\
\cA(\cC,-m-s/h)= \left\{\sum_{\substack{i\text{ an}\\\text{exponent}}}
    p_i\om_{-i}\mid p_i\in\C[t^{-1}],\quad \deg p_s=m,\quad \deg
    p_i\le\begin{cases} m &
    \text{for $i<s$}\\ m-1 & \text{for $i>s$}
  \end{cases}
  \right\}.
\]
Moreover, 
\[
\cA(\cC,-m-s/h)_{\mathrm{hom}}=\bC^\times t^{-m}\omega_{-s}. 
\]

In types $A$, $B$, and $C$, one can take
  $\om_{-s}=\om_{-1}^s$.  Finally, in types $A$ and $C$, one can simply
  view $\cA(\cC,r/h)$ as appropriate degree $r$ polynomials in $\om_{-1}$. In type $A_{n-1}$ (resp. $C_n$), the set of formal types consists of polynomials in $\omega_{-1}$ with no terms in degrees divisible by $n$ (resp. polynomials having no even degree terms).

\subsection{Fundamental strata and toral connections}\label{ss:Strata} In this subsection, we review some key aspects from the theory of fundamental strata for
formal connections and toral connections due to Bremer and
Sage \cite{BremerSageisomonodromy, BremerSagemodspace, BSregular,
  BSminimal, Sageisaac}.

\subsubsection{Naive leading term}

Consider the formal connection
\[
\nabla=d+(X_{-r}t^{-r}+X_{1-r}t^{1-r}+\dots)\dtt,
\] 
where $r$ is a positive
integer, $X_i\in\fg$, and $X_{-r}\ne 0$.  If $X_{-r}$ is
non-nilpotent, then the slope of $\nabla$ is $r$.  However, if
$\nabla$ has non-integral slope, then $X_{-r}$ is nilpotent; in fact, this
leading term will be nilpotent in any trivialisation of the underlying
bundle.

\subsubsection{} 
A more useful notion of the leading term of a formal connection
is obtained by the geometric theory of fundamental strata. Let $\cB$ be the Bruhat--Tits building
associated to the loop group $G(\cK)$.  Each $x\in\cB$ gives rise to
an $\bR$-filtration $\{\fg_{x,r}\}_{r\in\bR}$ of $\fg(\cK)$ by
$\bC[\![t]\!]$-lattices called the \emph{Moy--Prasad filtration
  associated to $x$}.  This grading is $1$-periodic in the sense that
  $\fg_{x,r+1}=t\fg_{x,r}$.  If we set
  $\fg_{x,r+}=\bigcup_{s>r}\fg_{x,s}$, then the set of $r$ for which
  $\fg_{x,r+}\ne \fg_{x,r}$ is discrete.  If $x$ is in the standard
  apartment determined by the split maximal torus $T(\cK)$, then the
  filtration is determined by a grading
  \[
\fg[t,t^{-1}]=\bigoplus_{r\in \mathbb{R}} \fg_x(r).
\]

\subsubsection{} 
  A $G$-stratum is a triple $(x,s,\beta)$, where $x\in\cB$, $s\ge 0$,
  and $\beta$ is a functional on $\fg_{x,s}/\fg_{x,s+}$, the
$s$th piece of the associated graded Lie algebra for the filtration.  Let $G_x$ be the parahoric subgroup corresponding to the facet
containing $x$, and let $G_{x+}$ be its pro-unipotent radical.  The
quotient $G_x/G_{x+}$ is then a complex reductive group.  This group
acts on the space of functionals $(\fg_{x,s}/\fg_{x,s+})^*$, and we
say that the $G$-stratum is \emph{fundamental} if $\beta$ is a
semistable point of this representation.

\subsubsection{} It will be convenient to reformulate the semistability
  criterion.  The Killing form induces an isomorphism
\begin{equation}\label{e:duality}
(\fg_{x,s}/\fg_{x,s+})^*\cong \fg_{x,-s}/\fg_{x,-s+},
\end{equation}
and a
functional $\beta$ is
semistable if and only if the corresponding coset consists entirely of
non-nilpotent elements.  Now, suppose that $x$ is in the standard
apartment, so that
$\fg_{x,-s}/\fg_{x,-s+}\cong\fg_x(-s)$.  We let
$\beta^0\in\fg_x(-s)$ be the image of $\beta$ under the composition of these two isomorphisms and call it the 
 homogeneous representative of $\beta$. 
Then, the functional $\beta$ is semistable if and only if $\beta^0$ is non-nilpotent.

\subsubsection{} \label{sss:slope}
Let $\nabla$ be a formal connection.  Upon fixing a trivialisation
$\phi$ of the underlying bundle, the connection can be written
$\nabla=d+A\dtt$ for some $A\in\fg(\cK)$.  

\begin{defe} We say that $\nabla$
contains the positive depth stratum $(x,s,\beta)$ with respect to the
trivialisation $\phi$  if $A\in\fg_{x,-s}$
and the coset $A+\fg_{x,-s+}$ corresponds to $\beta$ under the
isomorphism \eqref{e:duality}.\footnote{The definition of stratum
  containment when $r=0$ may be found in ~\cite{BremerSagemodspace,
    BSminimal}.}  In this case, $\beta$ is called the leading term of
$\nabla$ with respect to the Moy--Prasad filtration determined by $x$.
\end{defe} 

It is shown in \cite[Theorem 4.10]{BSminimal} that any irregular singular formal
connection $\nabla$ contains a fundamental stratum of positive depth
$r$, where $r$ is the slope of $\nabla$.
Moreover, this is the minimum possible depth of a stratum contained in
$\nabla$, and any other stratum of depth $r$ contained in $\nabla$ is
fundamental.  In particular, the slope of a formal connection may equivalently be
defined as the depth of any fundamental stratum contained in it.  This
is a useful point of view as it is often easier to compute the slope
using fundamental strata rather than using the definition given earlier.

\subsubsection{Regular strata}

A fundamental stratum $(x,r,\beta)$ is called \emph{regular} if every coset representative
of $\beta$ in $\fg_{x,-r}/\fg_{x,-r+}$ is regular semisimple.  If $x$
is the standard apartment, this is equivalent to $\beta^0$ being
regular semisimple.  The connected centralisers of these elements are
conjugate maximal tori, and we refer to the stratum as $S$-regular if
$S$ is one of these connected centralisers. 

 It turns out that $S$-regular strata of
depth $r$ exist only under stringent conditions on $S$ and $r$~\cite{BSregular}.
First, $S$ must correspond to a regular class in the Weyl group under the
bijection between classes of maximal tori in the loop group and
conjugacy classes in the Weyl group.  Second, $e^{2\pi i r}$ must be
a regular eigenvalue of this conjugacy class, i.e., an eigenvalue
corresponding to an eigenvector which is not on any wall in the
reflection representation.

\subsubsection{Toral connections} 
Suppose that
$\nabla$ is a connection containing a regular stratum $(x,r,\beta)$,
and assume $x$ is in the standard apartment. Let
$S=Z_{G(\cK)}(\beta^0)^\circ$ and 
$\fs:=\Lie(S)$.  In this case, we say that $\nabla$ is an \emph{$S$-toral connection}.  The
basic property of such a connection is that $\nabla$ is $G(\cK)$-gauge equivalent to $d+A\dtt$ with
$A\in\fs$. In fact, one can find an explicit rational canonical
form for $\nabla$. 

\begin{defe}  The set $\cA(S,r)$ of \emph{$S$-formal types of depth $r$}
is defined to be the subset of  $\bigoplus_{0\le s\le
  r}\fs(-s)$ consisting of elements with regular semisimple component in $\fs(-r)$.
  \end{defe} 
  
  The set $\cA(S,r)$ is
endowed with a natural action of the relative affine Weyl group
$\Waff_S\cong W_S\ltimes S/S_0$ on $\cA(S,s)$; here, 
$W_S$ acts by conjugation, and $S/S_0$ acts by translations along
$\fs(0)$.

\begin{thm}\cite{BSregular} Let $\nabla$ be an
  $S$-toral connection of slope $r$.  Then there exists $g\in G_{x+}$
  such that $g\cdot\nabla=d+A\dtt$ for some formal type $A\in\cA(S,r)$.  Moreover, two formal
  types give rise to isomorphic
formal connections if and only if they are in the same
$\Waff_S$-orbit.
\end{thm}

\subsection{Coxeter toral connections and the proof of
  Theorem~\ref{t:CoxToral}} \label{ss:Proof}

\begin{defe}\label{d:CoxToral}  A \emph{Coxeter toral connection} (or $\cC$-toral connection)
  is a formal connection which is toral with respect to a maximal torus conjugate to $\cC$. 
\end{defe}

It is well known that the regular eigenvalues of the Coxeter class in
$W$ are precisely the primitive $h$th roots of unity.  It thus
follows that the slope of a Coxeter toral connection has the form
$r/h$, where $\gcd(r,h)=1$.  Thus, a Coxeter toral connection is a formal
Coxeter connection in the sense of Definition~\ref{d:formCox}.

By the discussion above, every Coxeter toral
connection of slope $r/h$ is gauge equivalent to a
  connection of the form $d+A\dtt$ with  $A\in \cA(\cC,r/h)$.  It remains to  consider the action of the relative affine Weyl group on $\cA(\cC,r/h)$.

  The lattice component of the action of $\Waff_{\cC}$ is trivial
  because $\fc(0)=0$.  Moreover, it was shown in \cite{KSCoxeter} that
  $W_\cC\cong\mu_h$ with the action on $\cA(\cC,r/h)$ determined by $\zeta$ acting on $\fc(i/h)$ by
  $\zeta^{-i}$. Thus, the formal type $A$ is determined up to
  multiplication by an element of $\mu_h$.

It is now clear that Theorem~\ref{t:CoxToral} is equivalent to the
following statement:

\begin{thmbis}{t:CoxToral} A connection $\nabla$  is a formal Coxeter connection if
  and only if it is Coxeter toral.
\end{thmbis}

Note that we have already proved one implication, namely, Coxeter toral connections have slope $r/h$, and are therefore formal Coxeter connections. The rest of this subsection is devoted to proving the converse. 

\subsubsection{}
First, we show that very few Moy--Prasad filtrations admit semistable points of depth $r/h$
with $\gcd(r,h)=1$. 
\begin{prop} If $(x,r/h,\beta)$ is fundamental, then $x$ must be the
  barycentre of an alcove.  
\end{prop}
\begin{proof}  By equivariance of Moy--Prasad filtrations, it suffices
  to assume that $x$ is in the fundamental apartment.  Given $J\subset
 \Delta$, let  $W_J$ be the Weyl subgroup generated by the simple
  reflections from $\Delta_J$, and let $\check{\rho}_J$ be the half-sum of the
  positive coroots of the root system generated by $\Delta_J$ (viewed
  as an element of the fundamental apartment).  Reeder and
  Yu have shown in \cite[Theorem 8.3]{RY} that if $(x,r/m,\gamma)$
  is fundamental
  with $\gcd(r,m)=1$, then there exists $J\subset \Delta$
  such that $m$ divides a degree of $W_J$ and $x$ is
  conjugate under the affine Weyl group to an element of the affine
  subspace $\check{\rho}_J/m+\Delta_J^\perp$.

  We now specialise to the case where $m=h$. If $J$ is a proper subset of $\Delta$, then all degrees of
  $W_J$ are smaller than $h$. Hence, $J=\Delta$. As $\Delta^\perp=0$, we conclude that
  $x$ is in the affine Weyl group orbit of $\check{\rho}/h$, which is the
  barycenter of the standard alcove.
\end{proof}

Next, let $I\subset G(\cK)$ be the standard Iwahori subgroup associated to
$B$ and  $x_I$  the barycenter of the fundamental alcove.   By
the above proposition and the equivariance of stratum containment, we
can restrict to formal Coxeter connections containing a fundamental
stratum based at $x_I$.  We will show that all fundamental strata based at
$x_I$ are in fact \emph{regular} with leading term centralised by a
Coxeter torus.  

The 
Moy--Prasad grading at $x_I$ coincides with the standard Iwahori grading of \eqref{eq:Iwahori}, i.e.,  
\[
\fg_{x_I}(j/h)=\fg_{I}(j/h),\qquad \forall j\in \Z. 
\]
Fundamental strata of
depth $r/h$ based at $x_I$ are determined by non-nilpotent
elements of \mbox{$\fg_{x_I}(m+j/h)$}.  We find these elements by exploiting
the relationship between the Iwahori grading and the Coxeter grading
on $\fg$ considered in \S \ref{ss:CoxeterGrading}. 

Choose $u$ such that $u^h=t$, and consider a ramified cover
$\cK_b/\cK$ such that $\check{\rho}(u)$ is an element of $T(\cK_b)$.  (If $G$ is
of adjoint type, one can take $b=h$; otherwise, one needs to pass to
$b=2h$.) One obtains
\begin{equation}\label{eq:Kostant}
\Ad(\check{\rho}(u))\fg_{x_I}(m+j/h)=u^{mh+j}\fg_j.
\end{equation}

Kostant has shown that the
elements of $\fg_i$ are either nilpotent or regular semisimple and that
they are regular semisimple precisely when they have a nonzero
projection onto every root space of height congruent to $j$ modulo $h$ \cite{Kostant59}.
Moreover, if $X,X'\in\fg_i$ are regular semisimple, then $X'$ is
$T$-conjugate to $\lambda X$ for some nonzero $\lambda$.  It follows
that the analogous statements hold for the Iwahori grading.  In
particular, fundamental strata based at $x_I$ of depth $r/h$ (for
$\gcd(r,h)=1$) are actually regular.  Moreover, since $\beta^0\in\fg_{x_I}(-1/h)$ is regular semisimple, its centraliser in $\fg(\cK)$ is a Cartan subalgebra which corresponds to a regular conjugacy class $[w]$ in $W$. As 
$e^{-2\pi i r/h}$ is a regular eigenvalue of $[w]$, we must have that $[w]$ is the Coxeter
class.  Thus, we have shown that every connection of slope $r/h$ is
Coxeter toral.  This completes the proof of Theorem~\ref{t:CoxToral}. \qed

\subsubsection{Remark} Note that if a formal connection contains a semistable
stratum $(x_I,r/h,\beta)$ with respect to some trivialisation, then
after a further constant gauge change by an element of $T$, one can
assume that $\beta^0\in\fc(-r/h)$.  In other words, we can take the Coxeter torus
centralising $\beta^0$ to be $\cC$.

\subsection{Alternate proof of Theorem \ref{thm:Jordan}} \label{ss:RatJor}
Our goal is to show that for $A\in \cA(\cC,r/h)$, the
  connection $\nabla=d+A\dtt$ has Jordan form as in
  Theorem~\ref{thm:Jordan}.  Write $A=\sum_{i=1}^r
  A_i$, where $A_i\in\fc(-i/h)$ and $A_r\ne 0$.  Recall that $u^h=t$ and that we take a
  ramified cover $\cK_b$ of $\cK$ large enough so that $\check{\rho}(u)\in
  T(\cK_b)$. Also, $\ft'$ was defined to be the centraliser of the regular semisimple element $N+E$ (see \S \ref{ss:CoxeterGrading}). 
  
   We have already seen in \eqref{eq:Kostant} that $\Ad(\check{\rho}(u))(A_i)=u^{-i}X_i$
  for some $X_i\in\ft'_{-i}$; moreover, $X_{r}$ is regular
  semisimple.  Applying a gauge change by $\check{\rho}(u)$, we obtain
\[
d+\left[\left(\sum_{i=1}^r hu^{-i} X_{i}\right) -\check{\rho}\right]\duu.
\]

  Note that $\check{\rho}$ is orthogonal to $\ft'$ with respect to the Killing
  form.  It follows that $\check{\rho}$ is a sum of nonzero root vectors for
  $\ft'$, so it is possible to choose $Y\in\ft'$ such that
  $[Y,X_{-r}]=\check{\rho}$.  Gauge transformation by $\exp(u^r Y)$ gives
\[
d+\left[\left(\sum_{i=1}^r hu^{-i} X_{i}\right) +uD\right]\duu,
\]
  where $D\in\fg(\C[u])$.  It is now a standard result that this
  connection is $G(\cK_b)$-gauge equivalent to
\begin{equation} 
\label{e:JF} d+\left(\sum_{i=1}^r hu^{-i}
  X_{i}\right)\duu.
\end{equation} 

Since $X_{i}$ is nonzero only for exponents, we see
  that this connection has the desired Jordan type as in Theorem \ref{thm:Jordan}, albeit with $T$ replaced by $T'=C_G(N+E)$.  Moreover, if we let $\zeta_h$ be a primitive $h$th root of
  $1$, then $\theta=\rho(\zeta_h)$ normalises $T'$, satisfies $\theta^b=1$, and acts by
  $\zeta_h^i$ on $\ft'_i$.  It follows that $\theta$ induces a Coxeter
  element in the Weyl group $N_G(T')/T'$ and that the $X_{i}$ is an eigenvector of this Coxeter element with eigenvalue $\zeta_h^{-i}$. This concludes the proof. 
  \qed
  
\subsubsection{Remark} 
   Every Jordan form as in \eqref{e:JF} comes from a rational canonical form. Indeed, if $X_i\in\ft'_{-i}$ and $X_r$ is regular
  semisimple, one obtains the desired
formal type by setting 
  \[
  A=\frac{1}{h}\Ad(\check{\rho}(u^{-1}))
  \Big(\sum_{i=1}^r hu^{-i}
  X_{i}\Big).
  \]

\subsection{An algorithm for computing the rational canonical form}\label{s:algo} It is possible to find a rational canonical form for an $S$-toral connection by an explicit algorithm given in \cite{BremerSagemodspace, BSregular}.\footnote{When $G=\GL_n$ and $S=T(\cK)$, we recover the formal simplification of non-resonant irregular singular differential equations given in \cite[\S 11]{Wasow}.}   For later use, we give a brief
sketch of this algorithm for $\cC$-toral connections. 

\subsubsection{} 
Consider the formal Coxeter connection $\nabla=d+(Y+Z)\dtt$ of slope $r/h$,
where $Y\in\fg_{x_I}(-r/h)$ is regular semisimple and
$Z\in\fg_{x_I,(-r/h)+}$.   We know from Theorem~\ref{t:CoxToral} that
$\nabla$ is formally isomorphic to a connection of the form $d+A\dtt$
with $A\in \cA(\cC,r/h)$ and that the formal type $A$ is unique up to scalar
multiplication by an element of $\mu_h$.   Note that the formal type is unique
if the leading term is specified.

\subsubsection{} 
We now show how to find a formal type of $\nabla$ by applying a gauge transformation by an element of the Iwahori subgroup $I$.  
We have already seen that, by applying a
constant gauge change, we can assume without loss of generality that
$Y=A_{-r}\in\fc(-r/h)$.  For all $j>0$, let $I_{j/h}$ denote the $(j/h)$th
subgroup of the Moy--Prasad filtration on $I$, so in particular,
$I_+=I_{1/h}$. Suppose inductively that we have constructed an element
of $g_{r-s-1} \in I_+$ such that the matrix of $g_{r-s-1}\cdot\nabla$
has the form 
\[
A_{-r}+\dots+A_{-1-s}+Z_{-s}+Z'.
\]
 Here, $A_j$ is an element of $\fc(j/h)$
which is zero if $j>0$. Moreover, $Z_{-s}\in\fg_{x_I}(-s/h)$, and
$Z'\in\fg_{x_I,(-s/h)+}$.  It is shown in
\cite{BremerSagemodspace,BSregular} that one can find
$h_{r-s}\in I_{(r-s)/h}$ such that upon setting
$g_{r-s}=h_{r-s} g_{r-s-1}$, the matrix of $g_{r-s}\cdot\nabla$ has
  the form 
  \[
  A_{-r}+\dots+A_{-1-s}+A_{-s}+Z'',
  \] where
  $A_{-s}$ is an element of $\fc(-s/h)$ which is $0$ if $s<0$.  The desired element of
  $I_+$ is now $g=\lim g_k$.

  \subsubsection{} Note that once one has transformed all terms up to degree $-s/h$
  into $\fc$, they remain unchanged under further steps in the
  algorithm.  In particular, to determine the formal type, it is only
  necessary to carry out the first $r-1$ steps of the algorithm for
  computing $g$.  Thus, the algorithm
  to determine the formal type  finishes in finitely many steps.

  \subsubsection{} For $s\ge 0$, one can describe the effect of $h_{r-s}$ on $Z_{-s}$
  very explicitly.   Let $\langle,\rangle$ denote the Killing form.  It is shown in \cite{BremerSagemo
    dspace,BSregular} that $A_{-s}$ is
  the projection of $Z_{-s}$ onto $\fc(-s/h)$ under the form
  $\Res\langle,\rangle\dtt$.  As an example, consider
  $\fg=\mathfrak{sl}_{n}$.  Given $s\ge 0$, let $0\le \bar{s}<n$ be
  the residue of $s$ modulo $n$.  It is easy to see that
  $Z_{-s}\mid_{t=1}$ is an element of the $(n-\bar{s})$th ``circulant
  diagonal'' of $\mathfrak{sl}_{n}(\bC)$ (i.e., has entries only in
  the $(n-\bar{s})$th superdiagonal and the $\bar{s}$th subdiagonal).
  Let $a$ be the sum of the entries of $Z_{-s}\mid_{t=1}$.
  If one takes $\om_{-1}=N+t^{-1}E$ to be the usual element of
  $\fc(-1/n)$, then $Z_{-s}$ projects onto $\frac{a}{n}\om_{-1}^s$.


\section{Applications} \label{s:applications} 
In this section, we explore several applications of our main theorem (Theorem \ref{t:mainIntro}) to Coxeter $G$-connections. 
We begin by discussing generalised Frenkel--Gross connections in \S \ref{ss:FG}, where we prove Theorem \ref{t:FG} describing their differential Galois group. 
Similarly, in \S \ref{ss:Airy} we study Airy connections and establish Theorem \ref{t:Airy}. 
Next, we discuss the differential Galois groups of framed Coxeter connections (\S \ref{ss:frame}) and $\cC$-framed Coxeter connections (\S \ref{ss:Cframed}). In \S \ref{ss:generalisedCFramed}, we define a generalised version of these connections where we allow nontrivial residue terms. 
These results together lead to \S \ref{ss:Galois}, where we prove Theorem \ref{t:Galois}, providing a constructive approach to inverse differential Galois theory for Coxeter connections. Finally, in \S \ref{s:alt}, we give a simpler approach which addresses almost all cases of Theorem \ref{t:Galois}.

Given a connection $\nabla$ on $\bGm$, we write $\nabla_0$ (resp. $\nabla_\infty$) for the formal connections obtained by restricting $\nabla$ to the formal disk $\cDt_0$ (resp. $\cDt_\infty$).

\subsection{Generalised Frenkel--Gross connections}\label{ss:FG} Recall that a generalised Frenkel--Gross connection is a (Coxeter) connection on $\bGm$ with an irregular singularity at $0$ of slope $1/h$ and 
 a regular singularity at $\infty$. 
These connections have played an important role in the geometric Langlands program  \cite{FG, HNY, Zhu, HJ, ChenLi}. If we assume the monodromy is unipotent, then we obtain the original Frenkel--Gross connection \cite{FG}. (This follows from physical rigidity as mentioned below.) For $G=\SL_n$, one recovers the generalised Kloosterman connections (or hypergeometric connections with an irregular singularity of slope $1/n$) studied extensively by Katz \cite{Katz87, Katz90}.

Let $\nabla$ be a generalised Frenkel--Gross connection, and let 
$\rho_\nabla: \pi_1(\bGm)=\mathbb{Z} \to G$ denote its monodromy representation.  

\begin{thm} 
\begin{enumerate} 
\item  $\rho_\nabla(1)$ is regular. 
\item $\nabla$ is cohomologically rigid. 
\item $\nabla$ is physically rigid.\footnote{For $G=\GL_n$, Katz proved that physical rigidity is
equivalent to cohomological rigidity \cite[\S 1]{Rigid}. For more
general $G$, there exist cohomologically rigid $G$-connections
that are not physically rigid; see, for instance \cite{NAM, NBM}.} 
\item $\nabla$ is  $G(\cK$)-equivalent to
$\eta:=d+ a(N+S+Et^{-1})\dtt$ 
for some choice of $a\in \bC^\times$ and $S\in \ft$. 
\end{enumerate} 
\end{thm}

\begin{proof} 
Let
 $j:\bGm\hookrightarrow \bP^1$ be the natural inclusion. As shown in \cite{FG, KSCoxeter}, we have
\[
\dim H^1(\bP^1, j_{!*}\ad_\nabla) = \rankop(G) - \dim Z_\fg(\rho_\nabla(1)).
\]
Since dimensions are nonnegative, this immediately implies that $\dim Z_\fg(\rho_\nabla(1)) \le \rankop(G)$. The opposite inequality always holds, so equality must hold, proving (1). 

Next, by the discussion above, 
$H^1(\bP^1, j_{!*}\ad_\nabla) = 0$. Thus, $\nabla$ is cohomologically rigid. 

Part (3) is proved in \cite{Lingfei, HJ}. 

For (4),  in view of the discussion in \S \ref{ss:Rational}, we can choose $a\in \bC^\times$ such that 
\[
\nabla_0 \simeq d+a(N+Et^{-1})\dtt\simeq \eta_0. 
\]
Next, choose $S\in \ft$ such that $\exp(-2\pi i aS)$ is the semisimple part of the monodromy $\nabla_\infty(1)$. By \cite[Theorem 1.2]{Steinberg}, a regular class is determined by its semisimple part. Thus, we conclude 
\[
\rho_\nabla(1) = \textrm{the regular class in $G$ with semisimple part $\exp(-2\pi i aS)$} = \rho_\eta(1).\footnote{Here, we have slightly abused notation and written $\rho_\nabla(1)$ for the corresponding conjugate class.}  
\]
It follows that $\nabla_\infty\simeq \eta_\infty$.  By physical rigidity, we have $\nabla\simeq \eta$. 
\end{proof}

\subsubsection{Proof of Theorem \ref{t:FG} on the differential Galois group of $\nabla$}  First, observe that the semisimple part of the monodromy of $\nabla$ is $\exp(-2\pi i S)$. Since $S$ is torsion-free, Lemma \ref{l:connected} implies that the geometric monodromy group is connected. By Theorem \ref{t:mainIntro}, $G_\nabla$ is a reductive subgroup of $G$ of maximal degree. It remains to show that for $\sigma\in \Aut(G, \cP)$, we have  $\sigma(\nabla) \simeq \nabla$ if and only if $\sigma(S)$ is $\widetilde{W}$-conjugate to $S$ (with the usual modification for type $B_3$). 

Suppose $\sigma(\nabla) \simeq \nabla$. This implies that the
monodromies $\rho_{\nabla}$ and $\rho_{\sigma(\nabla)}$  are
$G$-conjugate, and so their semisimple parts $\exp(-2\pi i S)$ and $\exp(-2\pi i \sigma(S))$ are also conjugate. 
By Lemma \ref{l:WeylConjugate}, $S$ is $\widetilde{W}$-conjugate to $\sigma(S)$.  

Conversely, suppose $S$ and $\sigma(S)$ are in the same $\widetilde{W}$ orbit. By Lemma \ref{l:WeylConjugate}, $\exp(-2\pi i S)$  and $\exp(-2\pi i S')$ are $W$-conjugate. 
As we have seen, the monodromy is regular. Moreover, the
conjugacy class of a regular element is determined by the class
of its semisimple part. Thus, we conclude that the monodromies of $\sigma(\nabla)$ and
$\nabla$ are $G$-conjugate. Hence, 
$\sigma(\nabla)_\infty \simeq \nabla_\infty$. As noted in \S \ref{ss:Rational},  $\Aut(G, \cP)$ acts trivially on $\cA(\cC,1/h)$; thus, the 
formal types of $\nabla_0$ and $\sigma(\nabla)_0$ coincide. It follows that $\sigma(\nabla)_0 \simeq \nabla_0$. Physical rigidity then implies $\sigma(\nabla) \simeq \nabla$.  
This concludes the proof. \qed

\subsection{Airy connections} \label{ss:Airy} Recall that  an Airy $G$-connection is a (Coxeter) connection on $\bA^1$ with an irregular singularity at $0$ of slope $1+1/h$. 
For general simple groups $G$, these connections were first studied in \cite{KSCoxeter} and have
since played an important role in the geometric Langlands program and
the Deligne--Simpson problem \cite{KML, JY, YiAiry, KLMNS}. The original Airy equation gives rise to an Airy $\SL_2$-connection.

\begin{thm} Let $\nabla$ be an Airy connection. 
\begin{enumerate} 
\item $\nabla$ is cohomologically rigid. 
\item $\nabla$ is physically rigid. 
\item $\nabla$ is equivalent to  $d + at^{-1}(N + X + t^{-1}E)\,\dtt$ for some  $a \in \bC^\times$ and $X \in \fb$.
\end{enumerate} 
\end{thm} 

\begin{proof} These statements are proved in 
\cite{KSCoxeter}, \cite{HJ}, and \cite{YiAiry} respectively. 
\end{proof}

\subsubsection{Proof of Theorem \ref{t:Airy} on the differential Galois group} 
 First, observe that the monodromy of $\nabla$ is trivial. Thus, the geometric monodromy group is connected. By Theorem \ref{t:mainIntro}, $G_\nabla$ is a reductive subgroup of $G$ of maximal degree. Let $\FT(\nabla_0)$ be the (unique) formal type of $\nabla_0$ whose leading term (in the standard Iwahori grading \eqref{eq:Iwahori})  is 
$at^{-1}(N +  t^{-1}E)$. It remains to show that for $\sigma\in \Aut(G, \cP)$, we have  $\sigma(\nabla) \simeq \nabla$ if and only if $\FT(\sigma(\nabla_0))=\FT(\nabla_0)$ (with the usual modification for type $B_3$).

By physical rigidity, $\sigma(\nabla) \simeq \nabla$ if
  and only if $\sigma(\nabla)_0\simeq\nabla_0$.    Since
  \[
\sigma(\nabla) = d + at^{-1}(N + \sigma(X) + t^{-1}E)\,\dtt, 
  \]
 the formal types $\FT(\nabla_0)$ and
  $\FT(\sigma(\nabla_0))$  have the same leading term. By Theorem~\ref{t:CoxToral}, 
  $\sigma(\nabla)_0\simeq\nabla_0$ if and only if
  $\FT(\sigma(\nabla_0))=\FT(\nabla_0)$.  The result now follows from
  Theorem~\ref{t:mainIntro}.\qed

\subsubsection{More explicit description of $G_\nabla$}\label{sss:explicitX}
If $X$ is fixed by  $\Aut(G, \cP)$ (with suitable modifications in type $B_3$) then $\nabla$ is fixed by $\Aut(G, \cP)$, and $G_\nabla$ is the smallest reductive subgroup of $G$ of maximal degree. We now explain how to compute $G_\nabla$ explicitly in some cases where $X$ is not fixed by $\Aut(G, \cP)$.

In \S \ref{s:algo}, we gave an algorithm for computing $\FT(\nabla_0)$; in particular, this allows us to decide whether $\FT(\nabla_0)=\FT(\sigma(\nabla_0))$. This algorithm is effectively computable, but difficult to implement without a computer. 
However, if we choose $X\in \fb$ homogenous in the Iwahori grading, then it is easy to compute $\FT(\nabla_0)$ by hand up to $\deg(X)$. As we will see in the proposition below,  this usually suffices to determine whether $\FT(\nabla_0)=\FT(\sigma(\nabla_0))$.

\subsubsection{} Recall the definition of $s_0$ given in \S \ref{s0}. 
Let 
\[
\cM(s_0):=\bigoplus_{\substack{\alpha \in \Phi^+\\|\alpha|=h-s_0} }\fg_\alpha \subset \fb.
\]
In what follows, we consider those Airy connections for which $X\in \cM(s_0)$. 

\begin{prop}\label{p:AiryIGT} 
 There exists a nonempty
  Zariski-open subset of $\cM(s_0)$ for which the corresponding
Airy connections have $G_\nabla=G$. Moreover, in
type $D_4$, there is a nonempty subset of  $\cM(s_0)$ for which $G_\nabla$ has type $B_3$.  Specifically:
 \begin{enumerate} 
  \item Suppose we are in types $A_{2n-1}$ for $n\geq 2$, $D_n$ for
    $n\geq 5$, or
    $E_6$,
    and $\sigma\in \Aut(G, \cP)-\{1\}$.  Then,
    $G_\nabla=G$ if the projection of $X$ onto the $-1$-eigenspace of
    $\sigma$ in $\fc(-s_0/h)$ is nonzero. (Except for $D_{2n}$,
      this condition is simply that the projection of $X$ onto $\fc(-s_0/h)$
      is nonzero.)
  \item Suppose we are in type $D_4$. Then $G_\nabla=G$ if, for each 
    involution $\tau\in \Aut(G, \cP)$, the projection of $X$ onto the $-1$-eigenspace of
    $\tau$ in $\fc(-1/2)$ is
    nonzero.  On the other hand, if $\sigma(X)=X$ for one of the
    pinned involutions and, for another pinned involution $\tau$, the projection of $X$ onto the $-1$-eigenspace of $\tau$ in $\fc(-1/2)$ is
    nonzero, then $G_\nabla$ has type $B_3$. 
   \item If $G$ has type $B_3$, then $G_\nabla=G$  if
     the projection of $X$ onto the $-1$-eigenspace of $\tau$ in $\fc(-1/2)$ is
    nonzero for some pinned involution $\tau$ of $D_4\supset B_3$.
 \end{enumerate}
\end{prop}

\begin{proof} 
To determine $G_\nabla$, we will compute the
formal type of $\nabla_0$ up to degree $-s_0/h$.    As discussed in \S \ref{s:algo}, there
exists $g\in I_{(r-s_0)/h}$ such that 
\[
g\cdot \nabla_0=(a(N+t^{-1}E)t^{-1}+A_{s_0}+B)\dtt,
\]
 where
$A_{s_0}\in\fc(-s_0/h)$ and $B\in\fg_{x_I,-s_0/h+}$.  It follows that
\[
\FT(\nabla_0)= a(N+t^{-1}E)t^{-1}+A_{s_0}+\textrm{higher order terms}.
\]

Let $\sigma\in \Aut(G, \cP)$. Then
\[
\sigma(\nabla_0)=d+a(N+\sigma(X)+t^{-1}E)t^{-1}\dtt.
\]
We observe
that 
\[
\sigma(g)\cdot\nabla_0=d
+(a(N+t^{-1}E)t^{-1}+\sigma(A_{s_0})+\sigma(B))\dtt.
\]
 Thus, 
\[
\FT(\sigma(\nabla_0))=a(N+t^{-1}E)t^{-1}+\sigma(A_{s_0})+\textrm{higher order terms}.
\]

 Now, suppose we are in the types where $ \Aut(G, \cP)$ is generated by
 the involution $\sigma$.  Recall that $A_{s_0}$ is the projection of
 $X$ on $\fc(-s_0/h)$.  The condition on $X$ is precisely the
 condition that $\sigma(A_{s_0})\ne A_{s_0}$.  This implies that
 $\FT(\nabla_0)\ne\FT(\sigma(\nabla_0))$. By
 Theorem~\ref{t:Airy}, $G_\nabla=G$.

 The argument is similar in the remaining types and is left to the reader. 
\end{proof}

\subsubsection{} One can easily check the conditions of the above proposition explicitly. For example, in type $A$,  the condition is simply that the sum of the entries of $X$ is non-zero.

\subsection{Framed Coxeter connections}  \label{ss:frame} 
In this subsection, we give an explicit realisation of Coxeter connections 
 (Definition \ref{d:CoxConn}). To do so, we use the notion of framability,  which extends the theory of Fuchsian connections to allow irregular singularities \cite{BremerSagemodspace,Sagelocal}.

\begin{defe} \label{d:frame}
A connection $\nabla = d + M\, \dtt$ is called a \emph{framed Coxeter connection} on $\bGm$ if for some positive integer $r$ coprime to $h$, we have 
\[
M \in \fb \;\oplus\; \bigoplus_{j=1}^r \fg_I(-j/h) \subseteq \fg(\C[t^{-1}]),
\]
and the projection of $M$ onto $\fg_I(-r/h)$ is regular
semisimple.
\end{defe} 

Note that by a constant gauge change, one may assume that the leading term of $M$ lies in $\fc(-r/h)$.  
For instance, generalised Frenkel--Gross connections and Airy connections  are framable Coxeter connections
whose leading terms lie in $\fc(-1/h)$ and $\fc(-1-1/h)$, respectively.  

A framed Coxeter connection has a regular singularity at
$\infty$ because $M$ has no positive powers of $t$. As shown in \S
\ref{ss:Strata}, it is irregular at $0$ with slope $r/h$. Hence,
every framed Coxeter connection is a Coxeter connection in the sense of Definition
 \ref{d:CoxConn}.  We do not know if the converse holds in
general. However, as noted above, physical rigidity implies that the converse does hold when $r=1$ or
$r=h+1$.

\begin{prop} \label{p:Cox}
Let $\nabla$ be a Coxeter $G$-connection of slope $r/h$. 
\begin{enumerate} 
\item[(i)] If the geometric monodromy group $G_\nabla^{\mathrm{mono}}$
  is connected, then the differential Galois group $G_\nabla \subseteq
  G$ is a  reductive subgroup of maximal degree. Thus, if $G$ is of
  type $A_1$, $A_{2n}$, $B_n$ for $n\ge 4$, $C_n$, $E_7$, $E_8$, $F_4$, or $G_2$, then $G_\nabla = G$.
\item[(ii)] Let $\nabla=d+M\,\dtt$ be a framed Coxeter connection,  and assume the constant term $M_0$ of $M$ is torsion-free. Then the geometric monodromy group $G_\nabla^{\mathrm{mono}}$ is connected, and $G_\nabla\subseteq G$ is a reductive subgroup of maximal degree. In particular, this is the case if $M_0$ is nilpotent. 
\item[(iii)] The connection $\nabla$ is cohomologically rigid if and only if the dimension of the $G$-centraliser of $\rho_\nabla(1)$ equals $r \times \rankop(G)$.
\end{enumerate} 
\end{prop} 

\begin{proof} 
Part (i) follows from Theorem \ref{t:mainIntro}. Part (ii) follows from Lemma \ref{l:connected} and the fact that the semisimple part of the monodromy of $\nabla$ is $\exp(-2\pi i (M_0)_s)$. Part (iii) was proved in \cite{KSCoxeter}. 
\end{proof}

\subsubsection{More explicit description of $G_\nabla$} For the remainder of the subsection, we assume that we are in types $A_{2n-1}$, $D_n$,
$E_6$, or $B_3$. 
It is difficult to compute the differential Galois groups of
a framed Coxeter connection $\nabla=d+M\,\dtt$ in general.  However, as we now explain, one can put explicit
Zariski-open conditions on the second nonzero graded term (in the Iwahori grading) of $M$
guaranteeing that $G_\nabla=G$. This is a generalisation of the argument given in \S \ref{sss:explicitX}.

 Let $\hat{s}$ be a positive integer congruent modulo
$h$ to an exponent of $G$  that either does not appear,
or---in type $D_{2n}$---appears with smaller multiplicity, in any
proper reductive subgroup of maximal degree. Explicitly, in type
$A_{2n-1}$, $\hat{s}$ is congruent modulo $2n$ to an even number not
divisible by $n$; in type $D_{n}$,  $\hat{s}$ is congruent modulo
$2n-2$ to $n-1$; in type $E_6$, $\hat{s}$ is congruent modulo $12$ to
$4$ or $8$; and in type $B_3$,  $\hat{s}$ is congruent modulo $6$ to
$3$.  Note that the smallest possibility for $\hat{s}$ was denoted by $s_0$ in \S \ref{sss:explicitX}.

For convenience, we now restrict attention to framed Coxeter connections with 
\[
M\in \bigoplus_{j=1}^r \fg_I(-j/h).  
\] 
Assume $r\geq s_0$, and choose $\hat{s}<r$. 
 Consider
connections of the form
\begin{equation*} \nabla=d+(A_r+Y+Z)\,\dtt,
\end{equation*}
where $A_r$ is a nonzero element in
$\fc(-r/h)$, $Y$ is a nonzero element $\fg_{I}(-\hat{s}/h)$, and 
\[
Z\in\fg_{I} ((1-\hat{s})/h)\oplus\dots\oplus\fg_{I}(-1/h).
\]
 Let
$\cM(r,\hat{s})=\{(A_r,Y,Z)\}$ denote this parameter space.  For
$r>h$, we let
$\cM_{\mathrm{rf}}(r,\hat{s})$ denote the subset of triples
corresponding to residue-free framed Coxeter connections, i.e., framed
Coxeter connections on $\bP^1-\{0\}$.  Concretely, this means that $Y$ and $Z$ 
have no constant terms. Note that if we set $r=1+h$, $\hat{s}=s_0$, and $Z=0$, we obtain the set $\cM(s_0)$ of Proposition \ref{p:AiryIGT}.  

We have the following generalisation of Proposition \ref{p:AiryIGT}. 

\begin{prop}\label{p:framedIGT}  
    There exists a nonempty
  Zariski-open subset of $\cM(r,\hat{s})$ (resp. 
$\cM_{\mathrm{rf}}(r,\hat{s})$ if $r>h$) for which the corresponding
connections on $\bGm$ (resp. $\bP^1-\{0\}$) have $G_\nabla=G$.  In
type $D_4$, there is a nonempty subset of  $\cM(r,\hat{s})$ (resp. 
$\cM_{\mathrm{rf}}(r,\hat{s})$) for which $G_\nabla$ has type $B_3$.  Specifically:
   \begin{enumerate} 
  \item Suppose we are in types $A_{2n-1}$ for $n\geq 2$, $D_n$ for
    $n\geq 5$, or
    $E_6$,
    and $\sigma$ is the nontrivial pinned automorphism.  Then,
    $G_\nabla=G$ if the projection of $Y$ onto the $-1$-eigenspace of
    $\sigma$ in $\fc(-\hat{s}/h)$ is nonzero. (Except for $D_{2n}$,
      this condition is simply that the projection of $Y$ onto $\fc(-\hat{s}/h)$
      is nonzero.)
  \item If $G$ has type $D_4$, $G_\nabla$ is $G$ if for each pinned
    involution $\tau$, the projection of $Y$ onto the $-1$-eigenspace of
    $\tau$ in $\fc(-\hat{s}/h)$ is
    nonzero.  If $\sigma(Y)=Y$ and $\sigma(Z)=Z$ for one of the
    pinned involutions and, for another pinned involution $\tau$, the projection of $Y$ onto the $-1$-eigenspace of $\tau$ in $\fc(-\hat{s}/h)$ is
    nonzero, then $G_\nabla$ has type $B_3$. 
   \item If $G$ has type $B_3$, then $G_\nabla=G$ has type $B_3$ if
     the projection of $Y$ onto the $-1$-eigenspace of $\tau$ in $\fc(-\hat{s}/h)$ is
    nonzero for some pinned involution $\tau$ of $D_4\supset B_3$.
  \end{enumerate}
\end{prop}

\begin{proof} The proof is similar to that of Proposition \ref{p:AiryIGT}. First, 
  note that all connections determined by the parameter spaces
  $\cM(r,\hat{s})$ have nilpotent residue term, hence their
  differential Galois groups are reductive subgroups of maximal
  degree.
  
  We also observe that in each case (except when $G$ has type $D_4$ and $G_\nabla$ has type
  $B_3$), the condition on the parameter
  space  is that $Y$ lies in the complement of a proper linear
  subspace or the complement of the union of three proper linear
  subspaces.

To compute the differential Galois group, we will first compute the formal type at $0$ up to degree $-\hat{s}/h$.  As discussed in \S \ref{s:algo}, there
exists $g\in I_{(r-\hat{s})/h}$ for which  $g\cdot \nabla_0$
is in rational canonical form $d+(A_r+A_{\hat{s}}+B)\dtt$, where $A_{\hat{s}}\in\fc(-\hat{s}/h)$ and $B$
is the component of the formal type in degrees strictly greater than $-\hat{s}/h$.

Let $\sigma$ be a pinned automorphism. As noted in \S \ref{sss:action}, $\sigma(A_r)=
A_r$; thus, 
\[
\sigma(\nabla)=d+(A_r+\sigma(Y)+\sigma(Z))\dtt.
\]
 The
formal type of $\sigma(\nabla_0)$ is
$A_r+\sigma(A_{\hat{s}})+\sigma(B)$; indeed, it is obtained by applying the gauge change by
$\sigma(g)\in I_{(r-\hat{s})/h}$ to $\sigma(\nabla_0)$.

Suppose we are in cases where $ \Aut(G, \cP)$ is generated by
 the involution $\sigma$.  Recall that $A_{\hat{s}}$ is the projection of
 $Y$ on $\fc(-\hat{s}/h)$.  The condition on $Y$ in the statement of the proposition is precisely the
 condition that $\sigma(A_{\hat{s}})\ne A_{\hat{s}}$.  This implies that the
 formal types of $\nabla_0$ and $\sigma(\nabla_0)$ are not scalar multiples
 of each other. By Theorem \ref{t:CoxToral}, $\nabla_0$ is not (formally) gauge equivalent to $\sigma(\nabla_0)$. It follows that $\nabla$ is not gauge equivalent to $\sigma(\nabla)$. By Theorem \ref{t:mainIntro}, $G_\nabla=G$.

 The argument is similar in the remaining types, and is left to the reader. 
\end{proof}

\subsubsection{} 
The conditions of the above proposition are easy to check explicitly. For example, in type $A$,  the condition is simply that the sum of the entries of $Y\vert_{t=1}$ is non-zero. 

  \subsubsection{} A similar computation of the differential Galois
  group is possible for connections of the form $d+(A_r+Y+Z+W)\,\dtt$,
  where $W\in\fb$ and $A_r$, $Y$, and $Z$ are as in the proposition.
  However, in this case, one must further assume that the constant
  term of $A_r+Y+Z+W$ is torsion-free.

\subsection{$\cC$-framed Coxeter connections}\label{ss:Cframed}

The rational canonical forms of
Theorem~\ref{t:CoxToral} may be viewed as examples of framed Coxeter connections.  In this subsection, we compute their differential Galois groups. 

\begin{defe}
A \emph{$\cC$-framed Coxeter connection} is a connection on $\bGm$ of
the form $\nabla = d + A\dtt$, where $A \in \cA(\cC,r/h)$. 
\end{defe}

Note that a $\cC$-framed Coxeter connection $\nabla=d+A\dtt$ is, in particular, a framed Coxeter
connection. By definition, $A$ is a formal type of the formal connection
$\nabla_0$. Moreover, the constant term of $A$ is in $\bar{\fu}$,
where $\bar{\fu}$ denotes the nilpotent radical of the Borel subgroup
opposite to $\fb$. (This uses the fact that $\fg$ has trivial centre,
and hence $\fc(0)=\{0\}$.) It follows that $\nabla$ has a regular
singularity at $\infty$ with unipotent monodromy. The case where the monodromy is trivial is obtained 
by restricting to the \emph{residue-free} formal types
\[
\cA_\mathrm{rf}(\cC,r/h)\subset \cA(\cC,r/h).
\]
These are simply formal types with zero component in
$\fc(-s/h)$ for $s<h$.  They are precisely the formal types that give
rise to $\cC$-framed Coxeter connections on $\bP^1-\{0\}$.

\subsubsection{Differential Galois group} We now state the main result of this subsection. 

\begin{thm}\label{t:Cframed}  Let
  $\nabla=d+A\dtt$ be a $\cC$-framed Coxeter connection on $\bGm$. 
Then $G_\nabla$ is given by a simplified version of Table~\eqref{ta:mainIntro}, in which the conditions $\sigma(\nabla) \simeq \nabla$ are replaced by the requirement $\sigma(A)=A$ (and similarly for $\tau$ in type $D_4$).
\end{thm}

\begin{proof}
We first claim that $\nabla_0$ and
$\sigma(\nabla_0)$ are isomorphic if and only if 
$A=\sigma(A)$. Indeed, by Theorem \ref{t:CoxToral}, $\nabla_0$
and $\sigma(\nabla_0)$ are  isomorphic if and only if
$A=\zeta\sigma(A)$ for some $h$th root of unity
$\zeta$. Since $\sigma$ fixes the leading term of
$A$, i.e., the nonzero component of $A$ in
$\fc(-r/h)$, $A$ and $\sigma(A)$
differ by a nonzero scalar multiple if and only if they
are equal. This proves the claim. 

Now, suppose $G$ is of type $A_{2n-1}$ for $n\geq 2$, $D_n$ for $n\geq 5$,
or  $E_6$, and let $\sigma\in \Aut(G,
\cP)-\{1\}$. If $\sigma(A)\neq A$, then by the above discussion $\nabla$ and
$\sigma(\nabla)$ are not formally gauge equivalent at $0$. It
follows that they are not globally gauge equivalent. By Theorem
\ref{t:mainIntro}, we conclude $G_\nabla=G$. On the other hand, it is
clear that if $\sigma(A)=A$ then $\sigma(\nabla)=\nabla$; thus,
$G_\nabla$ is the unique proper reductive subgroup of maximal degree.

Next, suppose that $G$ has type $D_4$. If $\sigma(A)=A$ for some pinned involution $\sigma$, then $G_\nabla\subseteq G^{\sigma,\circ}$, so $G_\nabla$ has type $B_3$ or $G_2$. If $A$ is fixed by a pinned automorphism of order $3$, then it is of type $G_2$. If not, then the other involutions in $\Aut(G,\cP)$ do not fix $A$; thus, $\nabla$ and $\tau(\nabla)$ are not formally gauge equivalent at $0$ for involutions $\tau\neq \sigma$. Applying Theorem \ref{t:mainIntro} again, we conclude that $G_\nabla$ has type $B_3$. On the other hand, if $A$ is not fixed by any pinned automorphism of order $2$, then none of the resulting connections are formally isomorphic to $\nabla$ at $0$, and hence $G_\nabla=G$.

The argument for type $B_3$ is similar and is left to the reader.
\end{proof}

\subsubsection{The homogenous case} Suppose $A\in \cA_{\mathrm{hom}}(\cC, r/h)$ is a homogenous formal type. 
We studied the corresponding ($\cC$-framed) connections in \cite{KSCoxeter}, where we determined which of them are cohomologically rigid. 
In this case, $\sigma(A)=A$ for all $\sigma\in \Aut(G, \cP)$. Thus, in all types,  $G_\nabla$ is the smallest reductive subgroup of maximal degree in $G$.

\subsection{Generalised $\cC$-framed Coxeter connections}\label{ss:generalisedCFramed}
As in the case of generalised Frenkel--Gross connections, one can
obtain connections with non-unipotent monodromy by adding a diagonal
residue term to $\cC$-framed Coxeter connections: 
\begin{equation} 
\nabla = d + (A+S)\,\dtt, \qquad S \in \ft,\quad A\in \cA(\cC, r/h).
\end{equation} 
We call these \emph{generalised $\cC$-framed Coxeter
  connections}.  
These connections have a regular singularity at $\infty$ with semisimple monodromy
$\rho_\nabla(1)_s=\exp(-2\pi i S)$. The formal type at $0$ is unchanged
by the addition of $S$.

We have the following partial classification of differential Galois groups of the connections. 
\begin{prop}\label{p:genCframed} Suppose $S$ is torsion-free. Then $G_\nabla\subseteq G$ is a reductive subgroup of maximal degree. Moreover, we have:  
  \begin{enumerate} 
\item If $G$ is of type $A_1$, $A_{2n}$, $B_n$ for $n\ge 4$, $C_n$, $E_7$, $E_8$, $F_4$, or $G_2$, then $G_\nabla = G$. 
  \item If $G$ is of type $A_{2n-1}$ for $n\ge 2$, $D_n$ for $n\ge 5$, or $E_6$, then
  \begin{enumerate} 
  \item  If $\sigma(A+S)=A+S$ for $\sigma \in \Aut(G, \cP)-\{1\}$, then $G_\nabla$ is of type $C_n$, $B_{n-1}$, or $F_4$.  
  \item If $\sigma(A)\neq A$ or $\sigma(S)$ is not $\widetilde{W}$-conjugate to $S$, then $G_\nabla = G$.
  \end{enumerate}

  \item If $G$ is of type $D_4$, then
\begin{enumerate}
\item $G_\nabla$ has type $G_2$ if $\sigma(A+S) =A+S$ for all  $\sigma \in \Aut(G, \cP)$.
\item $G_\nabla$ has type $B_3$ if $\sigma(A+S) =A+S$ for exactly one involution $\sigma \in \Aut(G, \cP)$.
\item If $\sigma(A) \neq A$  or $\sigma(S)$ is not $\widetilde{W}$-conjugate to $S$ for all $\sigma \in \Aut(G, \cP)-\{1\}$, then $G_\nabla = G$.
\end{enumerate}

\item Suppose $G$ is of type $B_3$, and let $G'$ be a group of type $D_4$ containing $G$. Then 
\begin{enumerate} 
\item If $\tau(A+S)=A+S$ for some order three element $\tau \in \Aut(G', \cP')$, then $G_\nabla$ has type $G_2$. 

\item If $\tau(A) \neq A$  or $\tau(S)$ is not $\widetilde{W}$-conjugate to $S$ for some order three element $\tau \in \Aut(G', \cP')$, then $G_\nabla = G$. 
\end{enumerate} 
\end{enumerate}
\end{prop} 

Note that we do not have a complete understanding of the case in which $\sigma(S)$ is $\widetilde{W}$-conjugate (but not equal) to $S$. 

\begin{proof} If $\sigma(S+A)=S+A$ for some $\sigma\in \Aut(G, \cP)$, then $G_\nabla\subseteq G^{\sigma}$. On other hand, if $\sigma(A)\neq A$, then $\sigma(\nabla)_0\not\simeq  \nabla_0$ while if $\sigma(S)$ is not $\widetilde{W}$-conjugate to $S$, then $\sigma(\nabla)_\infty \not\simeq \nabla_\infty$. In either case, $\sigma(\nabla)\not\simeq \nabla$. The proposition readily follows from these observations and Theorem \ref{t:mainIntro}. 
\end{proof}

\subsection{Proof of Theorem \ref{t:Galois} on inverse differential Galois theory}\label{ss:Galois}  Our goal is to construct Coxeter connections of slope $r/h$ on $\bGm$ (or on $\bA^1$ if $r>h$) with all possible reductive subgroups of maximal degree. To prove the result for $\bGm$, we use Proposition~\ref{p:genCframed} to realise all possible differential Galois groups using 
generalised $\cC$-framed Coxeter connections whose underlying formal
type is homogeneous. For $\bA^1$, we obtain all possible differential Galois groups using the 
framed Coxeter connections considered in Proposition \ref{p:framedIGT}.

Throughout this subsection, we assume $G$ has a proper reductive subgroup of maximal degree. Recall the integer $s_0$ defined in \S \ref{s0}, and 
 the set $\cF\subset\ft$ defined in \S \ref{sss:SetF}.

\begin{prop}  For every reductive subgroup $K \subseteq G$ of maximal degree, there exists $S \in \cF$ such that the corresponding generalised $\cC$-framed Coxeter connection on $\bGm$
\[
\nabla = d + (A_r +S)\,\dtt
\]
has differential Galois group $G_\nabla = K$ for any
$A_r\in\fc(-r/h)$.
\end{prop}

 \begin{proof}
We will apply Proposition~\ref{p:genCframed}, noting that $A_r$ is
 always fixed by $\Aut(G,\cP)$.  
 
 We begin with the case where $K$ is the minimal reductive subgroup of
 maximal degree in $G$.  Here, we need only choose $S\in\cF$
 torsion-free and
 fixed by $\Aut(G,\cP)$ (with the usual modification
 in type $B_3$).  

Next, suppose we are in types $A_{2n-1}$ for $n\geq 2$, $D_n$ for $n\geq 5$, or $E_6$. Let $\sigma \in \Aut(G,\cP)-\{1\}$. 
Take $S\in \mathcal{F}-\cF^\sigma$ torsion-free. As $\sigma(S)\in
\mathcal{F}$,  Corollary~\ref{c:SetF} implies 
that $S$ and $\sigma(S)$ are not $\widetilde{W}$-conjugate to
$S$. Proposition~\ref{p:genCframed} now implies that $G_\nabla = G$.

Now, suppose $G$ is of type $D_4$. Let 
\[
S\in \cF - \bigcup_{\sigma \in \Aut(G, \cP) - \{1\}} \cF^\sigma
\]  
be a torsion-free element.  As above, this implies that $\sigma(S)$ is not $\widetilde{W}$-conjugate to $S$ for all nontrivial $\sigma \in \Aut(G, \cP)$, and hence \(G_\nabla = G\). 
On the other hand, consider a torsion-free element 
\[
S \in \cF^\sigma-\cF^{\tau},
\]
where $\sigma$ and $\tau$ are two distinct involutions
in $\Aut(G, \cP)$.  It follows that $S$ is fixed by exactly one pinned
involution, so by Proposition~\ref{p:genCframed}, $G_\nabla$ has type $B_3$.

The case $B_3$ is proved similarly and is left to the reader.

\end{proof}

\subsubsection{Connections on $\bA^1$}
Now suppose $r>h$, and consider the connections
\[
\nabla=d+(A_r+Y+Z)\,\dtt,
\]
 where 
 \[
 A_r\in \fc(-r/h)-\{0\},\quad Y\in\fg_I(-s_0/h)\cap
t^{-1}\fb, \quad \textrm{and} \quad Z\in t^{-1}\fb\cap(\fg_{I}((1-s_0)/h\oplus\dots\oplus\fg_{I}(-1/h)).
\]
The ``upper triangular''  component of  $\fg_{I}(-s/h)$ for $0<s<h$
is homogenous of degree $-1$ in $t$, so
this is a Coxeter connection whose sole singularity is an irregular singularity
of slope $r/h$ at $0$.

If we set $Y=Z=0$, then $G_\nabla$ is the minimal reductive subgroup of
 maximal degree in $G$.  Next, 
 assume that $Y\ne 0$, so the connections are parameterised by the set
 $\cF_{\mathrm{rf}}(r,s_0)$ considered in \S\ref{ss:frame}.
 Proposition~\ref{p:framedIGT} gives a nonempty Zariski-open subset of
 $\cF_{\mathrm{rf}}(r,s_0)$ (only involving the $Y$-component) which
 guarantees that $G_\nabla=G$.  In type $D_4$, the same result gives a
 nonempty subset for which $G_\nabla$ has type $B_3$.  This concludes
 the proof of Theorem \ref{t:Galois}. \qed

\subsection{Alternate approach using $\cC$-framed Coxeter
  connections}\label{s:alt}

 It turns out that one can prove almost all cases of Theorem \ref{t:Galois} using $\cC$-framed Coxeter connections with at most two non-zero components. 
The only cases requiring the more general connections considered above are: 
\begin{enumerate}
\item  $\bGm$:  $r=1$ in all types and $r<s_0$ in type
  $D_n$; 
  \item 
   $\bA^1$: $r=h+1$ in all types and $h<r< h+s_0$ in type
  $D_n$.
\end{enumerate}

\subsubsection{} It is a consequence of Theorem~\ref{t:Cframed} that, 
except in the cases mentioned above, $\cC$-framed Coxeter
  connections generically have $G_\nabla=G$. We now make this more precise by giving the explicit Zariski-open subsets for which this property holds.  
  
  Given $r>s_0$ (resp. $r>s_0+h$), define
  nonempty Zariski-open subsets of $\cA(\cC,r/h)$
  (resp. $\cA_\mathrm{rf}(\cC,r/h)$) by
  \begin{align*} \cA(\cC,r/h)^\circ&:=\{A\in\cA(\cC,r/h)\mid
                                     \text{$\sigma(A)\ne A$ for all
                                     involutions
                                     $\sigma\in\Aut(G,\cP)$}\},\text{ and}\\
    \cA_\mathrm{rf}(\cC,r/h)^\circ&:=\cA(\cC,r/h)^\circ\cap
                                    \cA_\mathrm{rf}(\cC,r/h).
  \end{align*}
  (Without the above conditions on $r$, these sets would be empty.) 
  
  In type $D_4$, note that the complements $(\cA(\cC,r/h)^\circ)^c$
 (resp. $(\cA_\mathrm{rf}(\cC,r/h)^\circ)^c$) are the spaces of
  formal types giving rise to connections on $\bGm$
  (resp. $\bP^1-\{0\}$) with proper differential
  Galois groups.  We define nonempty Zariski-open subsets of
  these parameter spaces by
\begin{align*} \cA(\cC,r/h)^{B_3}&=\{A\in\cA(\cC,r/h)\mid
                                     \text{$\sigma(A)\ne A$ for
                                   exactly one
                                     involution
                                     $\sigma\in\Aut(G,\cP)$}\},\text{ and}\\
    \cA_\mathrm{rf}(\cC,r/h)^{B_3}&=\cA(\cC,r/h)^{B_3}\cap
                                    \cA_\mathrm{rf}(\cC,r/h).
  \end{align*}

  \begin{prop} \label{p:Cframedexplicit} Suppose that $G$ has proper reductive subgroups of
    maximal degree.
  \begin{enumerate}\item  If $r>s_0$ (resp. $r>h+s_0)$, the subset of
    $\cA(\cC,r/h)$ (resp. $\cA_\mathrm{rf}(\cC,r/h)$) giving rise to
    $\cC$-framed connections with $G_\nabla=G$ is the nonempty,
    Zariski open subset $\cA(\cC,r/h)^\circ$
 (resp. $\cA_\mathrm{rf}(\cC,r/h)^\circ$).
  \item Suppose we are in type $D_4$.  Under the same conditions on
    $r$, the subset of
    $\cA(\cC,r/h)$ (resp. $\cA_\mathrm{rf}(\cC,r/h)$) giving rise to
    $\cC$-framed connections with $G_\nabla$ of type $B_3$ is the nonempty,
    Zariski locally closed subset $\cA(\cC,r/h)^{B_3}$
 (resp. $\cA_\mathrm{rf}(\cC,r/h)^{B_3}$).
\end{enumerate}
\end{prop}

\begin{proof}   The result follows immediately from Theorem~~\ref{t:Cframed}.
\end{proof}

\subsubsection{}
We now show how to obtain all possible differential Galois groups
(under the restrictions of Proposition~\ref{p:Cframedexplicit})
using explicit $\cC$-framed Coxeter connections with at most two non-zero
components.

 Let us first consider the case $\bGm$ with $r>s_0$.  One obtains the smallest $G_\nabla$ by
taking $d+A\dtt$, where $A$ is a homogenous formal type of depth $r/h$.\footnote{If
$r/h>1$, then this connection is nonsingular at $\infty$; if one takes
$A\in\cA(\cC,r/h)$ with nonzero components in degrees $-r/h$ and $-1/h$,
one again gets the smallest $G_\nabla$, and the connection will have
regular unipotent monodromy at $\infty$.}

The other possibilities arise by taking $A=A_r+A_{s_0}\in\cA(\cC,r/h)$
with nonzero components in degrees $r/h$ and $s_0/h$.  Outside of type $D_4$, one need only take
$A_{s_0}\ne 0$ to get $G_\nabla=G$.  In type $D_4$ (where $s_0=3$), one gets $G_\nabla=G$ by taking
$A_{3}$ not fixed by any pinned involution, and one gets $G_\nabla$ of type
$B_3$ if it is fixed by a single pinned involution.

Now, consider the case $\bA^1$ with 
$r>h+s_0$.  One obtains the smallest $G_\nabla$ by taking $A$ to be a
homogenous formal type of depth $r/h$.  The other possibilities arise
by taking $A=A_r+A_{h+s_0}\in\cA(\cC,r/h)$ with nonzero components in
degrees $-r/h$ and $-1-s_0/h$.  Outside of type $D_4$, one need only
take $A_{h+s_0}\ne 0$ to get $G_\nabla=G$.  In type $D_4$ (where $s_0+h=9$), one gets $G_\nabla=G$ by taking
$A_{9}$ not fixed by any pinned involution, and one gets $G_\nabla$ of type
$B_3$ if it is fixed by a single pinned involution.


  \end{document}